\documentclass{amsart}
\usepackage{amssymb,amsmath,latexsym}
\usepackage[dvips]{graphicx}

\newtheorem{thm}{Theorem}[section]
\newtheorem{dfn}[thm]{Definition}
\newtheorem{cor}[thm]{Corollary}

\newtheorem{lemma}[thm]{Lemma}

\newcommand{\del}{\backslash}
\newcommand{\NStep}{\ensuremath{N}}
\newcommand{\EStep}{\ensuremath{E}}
\newcommand{\ICal}{\ensuremath{\mathcal{I}}}

\newfont{\menutt}{cmtt8}

\title{Multi-Path Matroids}
\date{\today}
\author[Joseph E.~Bonin]
       {Joseph E.~Bonin}
\address[Joseph E.~Bonin]
{Department of Mathematics\\ The George Washington University\\
Washington, D.C. 20052, USA} \email [Joseph E.~Bonin]
{jbonin@gwu.edu}
\author[Omer Gim\'enez]
{Omer Gim\'enez$^\text{\dag}$}\thanks{\dag 
Partially supported by Beca Fundaci\'o Cr\`edit Andorr\`a.}
\address[Omer Gim\'enez]
{Departament de Matem\`atica Aplicada II\\ Universitat
Polit\`ecnica de Catalunya\\ Jordi Girona 1--3, 08034\\ Barcelona,
Spain} \email [Omer Gim\'enez] {Omer.Gimenez@upc.es}

\begin{document}

\begin{abstract}
  We introduce the minor-closed, dual-closed class of multi-path
  matroids.  We give a polynomial-time algorithm for computing the
  Tutte polynomial of a multi-path matroid, we describe their basis
  activities, and we prove some basic structural properties.  Key
  elements of this work are two complementary perspectives we develop
  for these matroids: on the one hand, multi-path matroids are
  transversal matroids that have special types of presentations; on
  the other hand, the bases of multi-path matroids can be viewed as
  sets of lattice paths in certain planar diagrams.
\end{abstract}

\maketitle

\section{Introduction}\label{sec:intro}

In~\cite{lpm1} it is shown how to construct, from a pair $P,Q$ of
lattice paths that go from $(0,0)$ to $(m,r)$, a transversal matroid
$M[P,Q]$ whose bases correspond to the paths from $(0,0)$ to $(m,r)$
that remain in the region bounded by $P$ and $Q$.  The basic
enumerative and structural properties of these matroids, which are
called lattice path matroids, are developed in~\cite{lpm1,lpm2}.  This
paper introduces multi-path matroids, a generalization of lattice path
matroids that share many of their most important properties.

Section~\ref{sec:def} starts by reviewing the definition of lattice
path matroids as well as an alternative perspective on these matroids
that uses collections of incomparable intervals in a linear order.
This alternative perspective leads to the starting point for
multi-path matroids: the linear order is replaced by a cyclic
permutation.  In addition to defining and providing examples of
multi-path matroids, Section~\ref{sec:def} also defines basic concepts
that are used in the rest of the paper.

Section~\ref{sec:basics} shows that the dual and all minors of a
multi-path matroid are multi-path matroids (lattice path matroids have
the corresponding properties; transversal matroids do not).  Proving
these properties entails developing several alternative presentations
for multi-path matroids.  In particular, we show that the bases of a
multi-path matroid can be viewed as certain sets of lattice paths in a
diagram (such as that in Figure~\ref{whirlpaths} on
page~\pageref{whirlpaths}) that has fixed global bounding paths and
one or more pairs of starting and ending points.

The diagrams we develop in Section~\ref{sec:basics} are crucial tools
in the next two sections.  Section~\ref{sec:tutte} shows that the
Tutte polynomial of a multi-path matroid can be computed in polynomial
time in the size of the ground set.  This result stands in contrast to
the following result of~\cite{cpv}: for any fixed algebraic numbers
$x$ and $y$ with $(x-1)(y-1)\ne 1$, the problem of computing
$t(M;x,y)$ for an arbitrary transversal matroid $M$ is $\#$P-complete.
Our work on the Tutte polynomial is cast in the general framework of
what we call computation graphs, which allow us to apply the idea of
dynamic programming to this computation.

Section~\ref{sec:activities} shows that, as is true of lattice path
matroids, internal and external activities of bases of multi-path
matroids have relatively simple lattice-path interpretations.  We also
sketch a somewhat faster, although more complex, algorithm for
computing the Tutte polynomial of a multi-path matroid via basis
activities.

The final section addresses several structural properties of
multi-path matroids.  For instance, we show that multi-path matroids
that are not lattice path matroids have spanning circuits and we
make some comments about minimal presentations of multi-path matroids.

We close this introduction by recalling several key notions;
see~\cite{ox} for concepts of matroid theory not defined here.  A
\emph{set system} is a multiset $\mathcal{A}=(A_j: j\in J)$ of subsets
of a finite set $S$.  A \emph{transversal} of $\mathcal{A}$ is a set
$\{x_j:j\in J\}$ of $|J|$ distinct elements such that $x_j$ is in
$A_j$ for all $j$ in $J$. A \emph{partial transversal} of
$\mathcal{A}$ is a transversal of a set system of the form $(A_k:k\in
K)$ with $K$ a subset of $J$.  Edmonds and Fulkerson~\cite{edmonds}
showed that the partial transversals of a set system $\mathcal{A}$ are
the independent sets of a matroid on $S$.  This matroid
$M[\mathcal{A}]$ is a \emph{transversal matroid} and the set system
$\mathcal{A}$ is a \emph{presentation} of $M[\mathcal{A}]$.  For a
basis $B$ of $M[\mathcal{A}]$, a \emph{matching of $B$ with
  $\mathcal{A}$} is a function $\phi:B\longrightarrow \mathcal{A}$
such that
\begin{itemize}
\item[(1)] $b$ is in $\phi(b)$ for each $b$ in $B$ and 
\item[(2)] the number of elements of $B$ that $\phi$ maps to any set
  $X$ in $\mathcal{A}$ is at most the multiplicity of $X$ in
  $\mathcal{A}$.
\end{itemize}
This terminology is suggested by the interpretation of set systems as
bipartite graphs~\cite[Section 1.6]{ox}.  In this paper, $\mathcal{A}$
will typically be an \emph{antichain}, that is, no set in
$\mathcal{A}$ contains another set in $\mathcal{A}$.  Presentations of
transversal matroids are generally not unique. A presentation
$(A_1,A_2,\ldots,A_r)$ of the transversal matroid $M$ \emph{contains}
the presentation $(A'_1,A'_2,\ldots,A'_r)$ of $M$ if $A'_i\subseteq
A_i$ for all $i$ with $1\leq i\leq r$.  We let $[n]$ denote the set
$\{1,2,\ldots,n\}$.

\section{Basic Definitions}\label{sec:def}

We start by reviewing lattice path matroids and an alternative
perspective on these matroids. The majority of this section is devoted
to defining multi-path matroids, providing illustrations, and defining
notation and concepts that are used in the rest of the paper.  Lattice
path matroids were introduced in~\cite{lpm1}; special classes of
lattice path matroids had been studied earlier from other perspectives
(see Section 4 of~\cite{lpm2}).

A lattice path can be viewed geometrically as a path in the plane made
up of unit steps East and North, or, more formally, as a word in the
alphabet $\{E,N\}$, where $E$ denotes the East step $(1,0)$ and $N$
denotes the North step $(0,1)$.  When viewed as a word in the alphabet
$\{E,N\}$, a lattice path does not have fixed starting and ending
points.  Thus, one may identify different geometric lattice paths that
arise from the same word; whether we identify such paths will depend
on, and should be clear from, the context.

Fix lattice paths $P$ and $Q$ from $(0,0)$ to $(m,r)$ with $P$ never
going above $Q$.  Let $\mathcal{P}$ be the set of all lattice paths
from $(0,0)$ to $(m,r)$ that go neither above $Q$ nor below $P$.  For
$i$ with $1\leq i\leq r$, let $N_i$ be the set
$$N_i=\{j\,:\, \text{ step } j \text{ is the } i\text{-th North step
  of some path in } \mathcal{P}\}.$$
The matroid $M[P,Q]$ is the
transversal matroid on the ground set $[m+r]$ that has
$(N_1,N_2,\ldots,N_r)$ as a presentation.  Note that $M[P,Q]$ has rank
$r$ and nullity $m$. A \emph{lattice path matroid} is any matroid
isomorphic to such a matroid $M[P,Q]$.

\begin{figure}
\begin{center}
\includegraphics[width = 3.0truein]{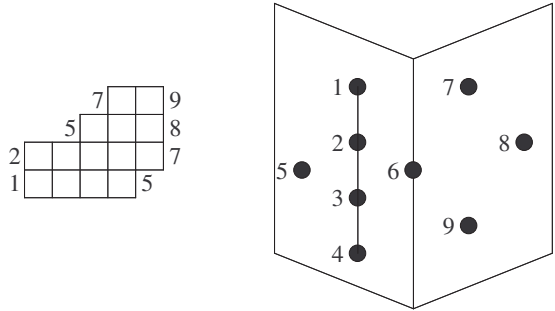}
\end{center}
\caption{A lattice path presentation and geometric representation
of a lattice path matroid.}\label{draw}
\end{figure}

Figure~\ref{draw} shows a lattice path matroid of rank $4$ and nullity
$5$.  The sets $N_1$, $N_2$, $N_3$, and $N_4$ are $\{1,2,3,4,5\}$,
$\{2,3,4,5,6,7\}$, $\{5,6,7,8\}$, and $\{7,8,9\}$. As this example
suggests, the sets $N_1, N_2,\ldots,N_r$ are intervals in $[m+r]$, and
both the left endpoints and the right endpoints form strictly
increasing sequences.  This motivates the following result
from~\cite{lpm2}.

\begin{thm}\label{thm:latticepath}
  A matroid is a lattice path matroid if and only if it is transversal
  and some presentation is an antichain of intervals in a linear order
  on the ground set.
\end{thm}

The following result~\cite[Theorem 3.3]{lpm1} starts to suggest a
deeper connection with lattice paths.

\begin{thm}\label{thm:bases}
  The map $R\mapsto \{i\,:\, \text{the } i\text{-th step of } R \text{
    is North}\}$ is a bijection from $\mathcal{P}$ onto the set of
  bases of $M[P,Q]$.
\end{thm}

Multi-path matroids are the generalizations of lattice path matroids
that result from using a cyclic permutation in place of the linear
order in Theorem~\ref{thm:latticepath}.

Fix a cyclic permutation $\sigma$ of the set $S$.  A
\emph{$\sigma$-interval} (or simply an \emph{interval}) in $S$ is a
nonempty subset $I$ of $S$ of the form
$\{f_I,\sigma(f_I),\sigma^2(f_I),\ldots,l_I\}$; this $\sigma$-interval
is denoted $[f_I,l_I]$ and the elements $f_I$ and $l_I$ are called,
respectively, the \emph{first} and \emph{last element of} $I$.  Note
that singleton subsets (which arise if $f_I$ is $l_I$) as well as the
entire set $S$ (which arises if $\sigma(l_I)$ is $f_I$) are
$\sigma$-intervals.  If one views the elements of $S$ placed around a
circle in the order given by $\sigma$, then the $\sigma$-intervals are
the sets of elements that can be covered by arcs of the circle; in the
case of a $\sigma$-interval $[f_I,l_I]$ that is $S$, the arc has a gap
between $l_I$ and $f_I$.

We now define our main object of study.

\begin{dfn}\label{def:multipath}
  A \emph{multi-path matroid} is a transversal matroid that has a
  presentation by an antichain of $\sigma$-intervals in some cyclic
  permutation $\sigma$ of the ground set.
\end{dfn}

The term ``multi-path'' comes from the alternative perspective on
these matroids given in Theorem~\ref{thm:altrep}.  To distinguish the
different types of presentations of interest in this paper,
presentations of the type in Theorem~\ref{thm:latticepath} are
\emph{interval presentations}, while those of the type in
Definition~\ref{def:multipath} are \emph{$\sigma$-interval
  presentations}.

Note that the first elements $f_{I_1},f_{I_2},\ldots,f_{I_r}$ that
arise from an antichain $\mathcal{I} = (I_1,I_2,\ldots,I_r)$ of
$\sigma$-intervals are distinct; thus, the rank of $M[\mathcal{I}]$ is
$r$, the number of intervals. Also, for $\mathcal{I}$ to be an
antichain of $\sigma$-intervals, the set $S$ can be in $\mathcal{I}$
only if $r$ is $1$.  However, Lemma~\ref{lem:relaxanti} shows that the
antichain condition in Definition~\ref{def:multipath} can be relaxed
without changing the resulting class of matroids; in some cases this
relaxation allows $S$ to be in $\mathcal{I}$.

In the following examples, $S$ is the set $[n]$ and $\sigma$ is the
cycle $(1,2,\ldots,n)$. Since a linear order can be ``wrapped around''
to obtain a cycle, every lattice path matroid is a multi-path matroid.
The converse is not true, as the first example shows.

\bigskip

\begin{figure}
\begin{center}
\includegraphics[width = 3.25truein]{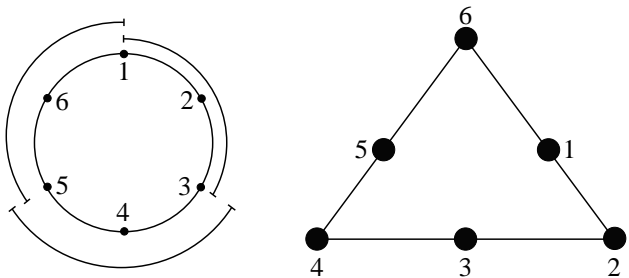}
\end{center}
\caption{The $3$-whirl $\mathcal{W}^3$ as a multi-path
matroid.}\label{whirlint}
\end{figure}

\noindent {\sc Example 1}.  The $3$-whirl $\mathcal{W}^3$ is an
excluded-minor for the class of lattice path matroids~\cite{lpm2};
however, Figure~\ref{whirlint} shows that the $3$-whirl is a
multi-path matroid.  The three intervals are $I_1=\{1,2,3\}$,
$I_2=\{3,4,5\}$, and $I_3=\{5,6,1\}$.  A similar construction shows that
all whirls are multi-path matroids.

\bigskip

There is only one interval presentation of a lattice path matroid
$M[P,Q]$ since $P$ and $Q$ correspond to, respectively, the greatest
and least bases in lexicographic order.  (See also~\cite[Theorem
5.6]{lpm2}.) In contrast, even lattice path matroids can have multiple
$\sigma$-interval presentations, as the next example shows.

\bigskip

\noindent {\sc Example 2}. All uniform matroids are lattice path
matroids. The following set systems are different $\sigma$-interval
presentations of the uniform matroid $U_{3,6}$ of rank $3$ on the set
$[6]$:
$$(\{1,2,3,4\},\{2,3,4,5\},\{3,4,5,6\}),$$
$$(\{1,2,3,4,5\},\{2,3,4,5,6\},\{3,4,5,6,1\}).$$

\bigskip

A presentation $\mathcal{A}$ of a transversal matroid $M$ is
\emph{minimal} if no other presentation of $M$ is contained in
$\mathcal{A}$.  Interval presentations of lattice path matroids are
minimal~\cite[Theorem 6.1]{lpm2}. We next show that multi-path
matroids that are not lattice path matroids can have multiple minimal
presentations that are $\sigma$-interval presentations.

\bigskip

\begin{figure}
\begin{center}
\includegraphics[width = 1.4truein]{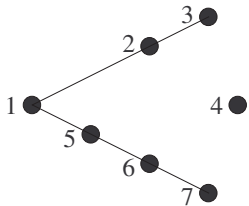}
\end{center}
\caption{A multi-path matroid that has multiple minimal
presentations.}\label{nonuni}
\end{figure}

\noindent {\sc Example 3}. The following set systems are
$\sigma$-interval presentations of the matroid shown in
Figure~\ref{nonuni} and both are minimal:
$$(\{5,6,7,1,2\},\{2,3,4\},\{4,5,6,7\}),$$
$$(\{6,7,1,2,3\},\{2,3,4\},\{4,5,6,7\}).$$

\bigskip

The rest of this section contains observations about multi-path
matroids as well as definitions that are used in later sections.

Many constructions, such as minors, involve subsets of the ground set;
in such settings, we use the following definition. The cyclic
permutation $\sigma$ on $S$ induces a cyclic permutation $\sigma_X$ on
each subset $X$ of $S$ defined as follows: for $x$ in $X$, the image
$\sigma_X(x)$ is the first element in the list
$\sigma(x),\sigma^2(x),\sigma^3(x),\ldots$ that is in $X$.  Thus,
$\sigma_X$ is formed from $\sigma$ by skipping over the elements that
are not in $X$.

There is an induced cycle on the $\sigma$-intervals in an antichain
$\mathcal{I}=(I_1,I_2,\ldots,I_r)$ of $\sigma$-intervals.  Indeed, the
last elements $l_{I_1},l_{I_2},\ldots,l_{I_r}$ are distinct since
$\mathcal{I}$ is an antichain, so a cyclic permutation $\Sigma$ on
$\mathcal{I}$ is given by $\Sigma(I_j) = I_k$ if $\sigma_X(l_{I_j}) =
l_{I_k}$ where $X$ is $\{l_{I_1},l_{I_2},\ldots,l_{I_r}\}$.  Likewise
the cycle $\sigma_Y$ on $Y=\{f_{I_1},f_{I_2},\ldots,f_{I_r}\}$ induces
a cyclic permutation $\Sigma'$ on $\mathcal{I}$.  The assumption that
$\mathcal{I}$ is an antichain gives the equality $\Sigma=\Sigma'$. We
use $\Sigma$ to denote the cyclic permutation of $\mathcal{I}$ induced
in this manner from $\sigma$.  For instance, in Example 1, $\Sigma$ is
$(I_1,I_2,I_3)$.

Fix an element $x$ in a $\sigma$-interval $I$.  It will be useful to
consider the two parts in which $I-x$ naturally comes. The \emph{first
  part of} $I-x$ is the empty set if $x$ is $f_I$, or the
$\sigma$-interval $[f_I, \sigma^{-1}(x)]$ if $x$ is not $f_I$.
Similarly, the \emph{last part of} $I-x$ is the empty set if $x$ is
$l_I$, or the $\sigma$-interval $[\sigma(x),l_I]$ if $x$ is not $l_I$.
From the set $I-x$ alone, references to the first and last parts could
be ambiguous (for instance, if $x$ is $f_I$ or $l_I$), but $x$ will be
clear from the context, so no confusion should result.

Note that an element $x$ in $S$ is a loop of $M[\mathcal{I}]$ if and
only if $x$ is in no interval in $\mathcal{I}$.  Thus if $x$ is a
loop, then the intervals in $\mathcal{I}$ are intervals in the linear
order $x<\sigma(x)<\sigma^2(x)<\cdots<\sigma^{-1}(x)$, so
$M[\mathcal{I}]$ is a lattice path matroid.  Likewise, note that if
some first element $f_I$ is not in $\Sigma^{-1}(I)$, then
$M[\mathcal{I}]$ is a lattice path matroid.

\section{Minors, Duals, and the Lattice Path
  Interpretation}\label{sec:basics}

This section shows that the class of multi-path matroids is closed
under minors and duals.  (Analogous properties hold for lattice path
matroids but not for arbitrary transversal matroids.)  We also develop
several alternative descriptions of multi-path matroids, some of which
involve lattice paths and so account for the name.  The lattice path
interpretations as well as closure under contractions enter into the
proof of closure under duality.

We start with a simple lemma that applies to all transversal matroids.

\begin{lemma}\label{lem:trimpres}
  Assume $X$ and $Y$ are in a set system $\mathcal{A}$ with
  $X\subseteq Y$.  Let $z$ be in $X$ and let $\mathcal{A}'$ be
  obtained from $\mathcal{A}$ by replacing one or more occurrences of
  $Y$ by $Y-z$.  Then $M[\mathcal{A}]= M[\mathcal{A}']$.
\end{lemma}

\begin{proof}
  Note that it suffices to prove the result in the case that one
  occurrence of $Y$ is replaced by $Y-z$, and for this it suffices to
  show that for any basis $B$ of $M[\mathcal{A}]$ and matching
  $\phi:B\longrightarrow\mathcal{A}$, we can find a matching
  $\phi':B\longrightarrow\mathcal{A}'$.  Clearly there is such a
  matching $\phi'$ if $z$ is not in $B$, or if $z$ is in $B$ but
  $\phi(z)$ is not $Y$.  Thus assume that $z$ is in $B$ and $\phi(z)$
  is $Y$. If $X$ is not in the image of $\phi$, then the map $\phi'$
  that agrees with $\phi$ except that $\phi'(z)$ is $X$ is the
  required matching.  Now assume $\phi(x)$ is $X$ for some $x$ in $B$.
  Since $x$ is in $X$ and therefore in $Y$, the following map $\phi'$
  is the required matching:
\begin{eqnarray}
 & &\phi'(w) =
   \left\{
     \begin{array}{ll}
      X, &\mbox{if $w=z$;}\\
      Y-z, &\mbox{if $w=x$;}\\
      \phi(w), &\mbox{otherwise.}
    \end{array}
   \right.  \notag
\end{eqnarray}
\end{proof}

It is well known and easy to see that if
$\mathcal{A}=(A_1,A_2,\ldots,A_r)$ is a presentation of a transversal
matroid $M$ on $S$, then any single-element deletion $M\del x$ is
transversal and $\mathcal{A}'=(A_1-x,A_2-x,\ldots,A_r-x)$ is a
presentation of $M\del x$.  Since deleting $\emptyset$ from any set
system in which it appears does not change the associated transversal
matroid, we may assume that $\emptyset$ is not in $\mathcal{A}'$.
Note that if $\mathcal{A}$ is an antichain of $\sigma$-intervals, then
the sets in $\mathcal{A}'$ are $\sigma_{S-x}$-intervals, but there may
be containments among these sets.  This issue is addressed through the
next lemma, which gives a relaxation of the antichain criterion in
Definition~\ref{def:multipath}.

\begin{lemma}\label{lem:relaxanti}
  Assume the transversal matroid $M$ has a presentation by a multiset
  $\mathcal{A}$ of $\sigma$-intervals that satisfies the following
  condition:
  \begin{quote}
   \textbf{(C)} \ if $I\subseteq J$ for $I,J\in\mathcal{A}$, then
   either $f_J$ or $l_J$ is in $I$.
  \end{quote}
  Then $M$ is a multi-path matroid and $\mathcal{A}$ contains a
  $\sigma$-interval presentation of $M$.
\end{lemma}

\begin{proof}
  If $\mathcal{A}$ is an antichain, there is nothing to prove, so
  assume $I$ and $J$ are in $\mathcal{A}$ and $I\subseteq J$.  By
  condition (C) and symmetry, we may assume $f_J$ is in $I$. By
  replacing $J$ if needed, we may assume no $\sigma$-interval in
  $\mathcal{A}$ whose first element is $f_J$ properly contains $J$.
  Let $\mathcal{A}'$ be the set system obtained from $\mathcal{A}$ by
  replacing $J$ by the $\sigma$-interval $J-f_J$, or eliminating $J$
  if $J-f_J$ is empty.  Lemma~\ref{lem:trimpres} gives the equality
  $M[\mathcal{A}]=M[\mathcal{A}']$; we will show that $\mathcal{A}'$
  satisfies condition (C).  The presentation of $M$ by
  $\sigma$-intervals that results from applying this modification as
  many times as possible must be an antichain, which proves the lemma.

  To show that $\mathcal{A}'$ satisfies condition (C), first note that
  the only pairs of intervals that potentially could contradict
  condition (C) must include $J-f_J$. Let $K$ be another interval in
  $\mathcal{A}'$.  If the containment $K\subseteq J-f_J$ holds, then
  $K$ is a subset of $J$ but does not contain $f_J$; it follows from
  condition (C) applied to $J$ and $K$ in $\mathcal{A}$ that $l_J$
  (which is also $l_{J-f_J})$ must be in $K$, as needed.  Now assume
  the containment $J-f_J\subseteq K$ holds.  If $l_K$ is in $J-f_J$,
  there is nothing to show, so assume this is not the case.  Since $J$
  is the largest set in $\mathcal{A}$ that has $f_J$ as its first
  element, $f_K$ is not $f_J$.  If $\sigma^{-1}(f_J)$ were in $K$,
  then $J$ and $K$ would contradict condition (C) for $\mathcal{A}$.
  Thus, the first element of $K$ must be $\sigma(f_J)$, so $f_K$ is in
  $J-f_J$, as needed.
\end{proof}

It is easy to check that if $(I_1,I_2,\ldots,I_r)$ is an antichain of
$\sigma$-intervals, then the set system $(I_1-x,I_2-x,\ldots,I_r-x)$
satisfies condition (C) of Lemma~\ref{lem:relaxanti}.  This
observation along with the remarks before that lemma prove the
following theorem.

\begin{thm}\label{thm:del}
  The class of multi-path matroids is closed under deletion.
\end{thm}

To show that the class of multi-path matroids is closed under
contractions, we give presentations of single-element contractions
(Lemma~\ref{lem:con}) that we then show satisfy condition (C) of
Lemma~\ref{lem:relaxanti}.

\begin{lemma}\label{lem:con}
  Let the antichain $\mathcal{I}$ of $\sigma$-intervals be a
  presentation of $M$, and let $\Sigma$ be the cycle
  $(I_1,\ldots,I_t,I_{t+1},\ldots,I_r)$ where $I_1,\ldots,I_t$ are the
  $\sigma$-intervals that contain a given element $x$.  A presentation
  of the contraction $M/x$ is given by:
\begin{itemize}
 \item[(a)]
 $\mathcal{I}$, for $t=0$;
 \item[(b)]
 $(I_2,I_3,\ldots,I_r)$, for $t=1$;
 \item[(c)]
 $\mathcal{I}':= \bigl((I_1\cup I_2)-x,(I_2\cup I_3)-x,\ldots,
 (I_{t-1}\cup I_t)-x,I_{t+1},\ldots,I_r\bigr)$, for $t>1$.
\end{itemize}
\end{lemma}

\begin{proof}
  Part (a) holds since $x$ is a loop of $M$ if $t$ is $0$. If $t$ is
  positive, then $x$ is not a loop, so the bases of $M/x$ are the
  subsets $B$ of $S-x$ such that $B\cup x$ is a basis of $M$. Part (b)
  follows since matchings $\phi:B\cup x\longrightarrow \mathcal{I}$
  map $x$ to $I_1$. For part (c), we need to show that for subsets $B$
  of $S-x$, there is a matching $\phi:B\cup
  x\longrightarrow\mathcal{I}$ if and only if there is a matching
  $\phi':B\longrightarrow\mathcal{I}'$.

  Assume first that $\phi:B\cup x\longrightarrow\mathcal{I}$ is a
  matching.  Assume $\phi(x)$ is $I_h$ and $\phi(b_i)$ is $I_i$ for
  all $i$ with $i\ne h$ and $1\leq i\leq r$.  The necessary matching
  $\phi'$ is given by
  \begin{eqnarray}
   & &\phi'(b_i) =
   \left\{
     \begin{array}{ll}
      (I_i\cup I_{i+1})-x, &\mbox{if $1\leq i<h$;}\\
      (I_{i-1}\cup I_i)-x, &\mbox{if $h< i\leq t$;}\\
      I_i, &\mbox{if $t<i\leq r$.}
    \end{array}
   \right.  \notag
  \end{eqnarray}

  Now assume $\phi':B\longrightarrow\mathcal{I}'$ is a matching and
  $\phi'(b_i)$ is $(I_i\cup I_{i+1})-x$ for $1\leq i <t$.  Since $x$
  is in $I_1,I_2,\ldots,I_t$, to complete the proof it suffices to
  construct an injection $\psi:\{b_1,b_2,\ldots,b_{t-1}\}
  \longrightarrow \{I_1,I_2,\ldots,I_t\}$ with each $b_i$ in
  $\psi(b_i)$.  Toward this end, classify $b_1,b_2,\ldots,b_{t-1}$ as
  follows: $b_i$ is a {\em leader} if it is in the first part of
  $I_i-x$, otherwise $b_i$ is a {\em trailer}. Note that if $b_i$ is a
  leader, then $b_i$ is in the first part of $I_j-x$ for every $j$
  with $1\leq j \leq i$.  Similarly, if $b_i$ is a trailer, then $b_i$
  is in the last part of $I_j-x$ for every $j$ with $i+1\leq j \leq
  t$. Define $\psi$ as follows: scan $b_1,b_2,\ldots,b_{t-1}$ in this
  order and for each leader $b_i$, let $\psi(b_i)$ be the first set
  among $I_1,I_2,\ldots,I_t$ that is not already in the image of
  $\psi$; then scan $b_{t-1},b_{t-2},\ldots,b_1$ in this order and for
  each trailer $b_i$, let $\psi(b_i)$ be the last set among
  $I_1,I_2,\ldots,I_t$ not already in the image of $\psi$.  Clearly
  $\psi$ is injective and $b_i$ is in $\psi(b_i)$ for all $i$.
\end{proof}

With this lemma, we can now complete our work on contractions.

\begin{thm}\label{thm:con}
  The class of multi-path matroids is closed under contraction.
\end{thm}

\begin{proof}
  We use the notation of Lemma~\ref{lem:con}.  It suffices to show
  that $M/x$ is a multi-path matroid.  This follows easily from parts
  (a) and (b) of Lemma~\ref{lem:con} if $t$ is at most $1$, so assume
  $t$ exceeds $1$.  To show that $M/x$ is a multi-path matroid, it
  suffices to show that $\mathcal{I}'$ satisfies condition (C) of
  Lemma~\ref{lem:relaxanti}.  To consider the sets in $\mathcal{I}'$
  as $\sigma_{S-x}$-intervals, we need only specify the endpoints of
  any set that is $S-x$.  If the $\sigma$-interval $I_h$ is $S-x$,
  where $t< h\leq r$, we take $I_h$ to be the $\sigma_{S-x}$-interval
  $[\sigma(x),\sigma^{-1}(x)]$.  If $(I_i\cup I_{i+1})-x$ is $S-x$, we
  take this to be the $\sigma_{S-x}$-interval
  $[f_{I_i},\sigma^{-1}(f_{I_i})]$.  Note that there are only three
  possible containments among the sets in $\mathcal{I}'$:
\begin{itemize}
\item[(i)] $(I_i\cup I_{i+1})-x\subseteq (I_j\cup I_{j+1})-x$ with
  $1\leq i,j<t$,
\item[(ii)] $(I_i\cup I_{i+1})-x\subseteq I_h$ with
  $1\leq i<t$ and $t< h\leq r$, and
\item[(iii)] $I_h\subseteq (I_i\cup I_{i+1})-x$ with
  $1\leq i<t$ and $t< h\leq r$.
\end{itemize}
In case (i), note that if $f_{I_j}$ is not in $(I_i\cup I_{i+1})-x$,
then $j<i$.  It follows that $l_{I_{j+1}}$ is in the $\sigma$-interval
$[\sigma(x),\sigma^{-1}(l_{I_{i+1}})]$, so $l_{I_{i+1}}$ is not in
$(I_j\cup I_{j+1})-x$.  This contradicts the assumed containment, so
$f_{I_j}$ is in $(I_i\cup I_{i+1})-x$ and condition (C) holds in case
(i).  Note that the containment in case (ii) holds only if $I_h$ is
$[\sigma(x),\sigma^{-1}(x)]$, so condition (C) clearly holds in this
case also.  Lastly, consider the containment in case (iii).  Since $x$
is not in $I_h$, if $f_{I_i}$ were not in $I_h$, then $I_h$ would be
either contained in or disjoint from
$[\sigma(f_{I_i}),\sigma^{-1}(x)]$, so either $I_h\subseteq I_i$ or
$I_h\subseteq I_{i+1}$ would hold.  That both conclusions are contrary
to $\mathcal{I}$ being an antichain shows that $f_{I_i}$ is in $I_h$,
so condition (C) of Lemma~\ref{lem:relaxanti} holds.  Thus, $M/x$ is a
multi-path matroid.
\end{proof}

\begin{figure}
\begin{center}
\includegraphics[width = 2.3truein]{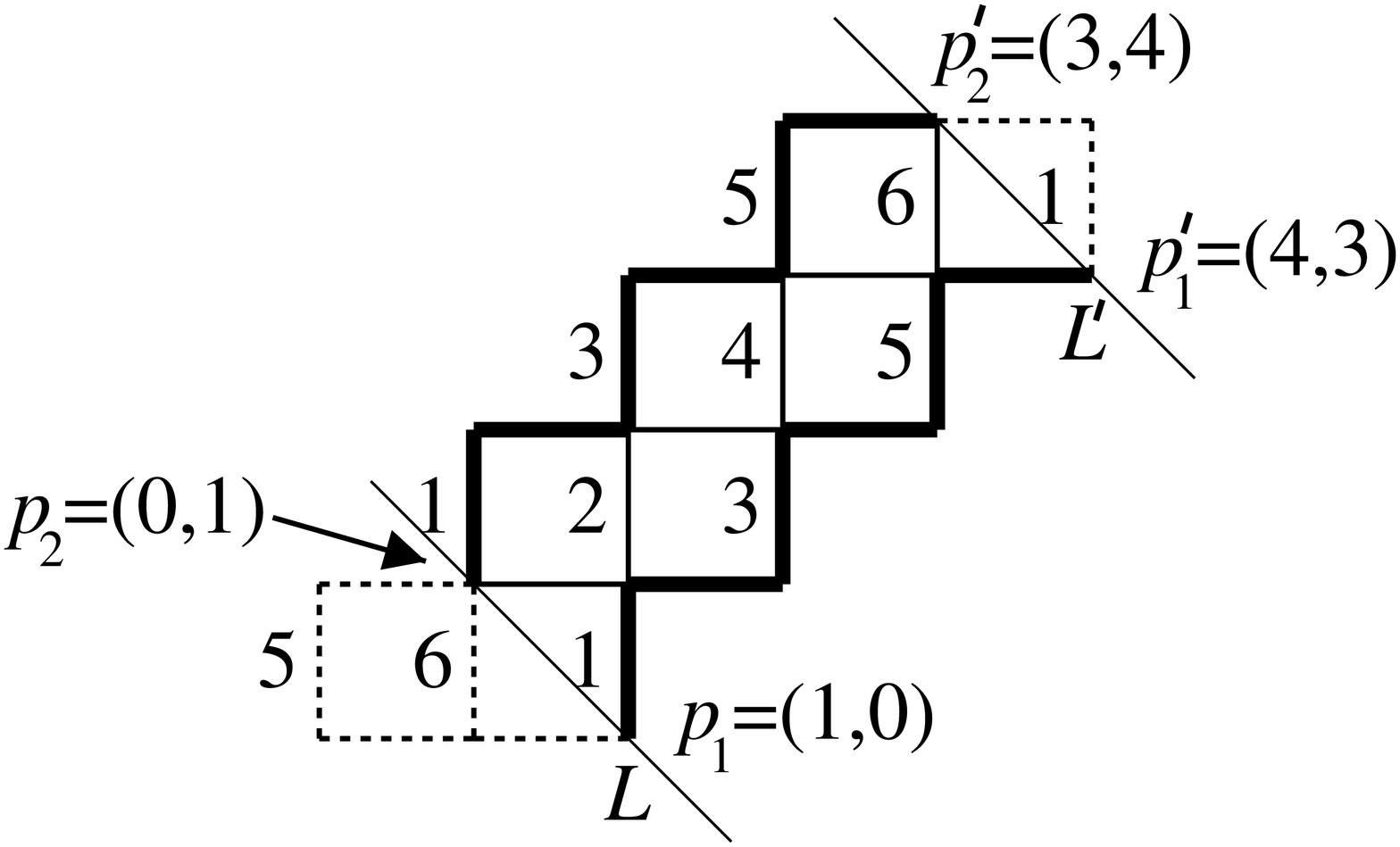}
\end{center}
\caption{The diagram $D(\mathcal{W}^3,1)$ with the labels on the
  North steps.}\label{whirlpaths}
\end{figure}

We now give an alternative perspective on multi-path matroids that
accounts for the name, extends the path interpretation of lattice path
matroids, and plays a pivotal role in much of the rest of this paper.
Figure~\ref{whirlpaths} illustrates these ideas with the $3$-whirl
(Example 1 in Section~\ref{sec:def}).  Assume $M[\mathcal{I}]$ has
rank $r$ and nullity $m$.  Fix an element $x$ of $M[\mathcal{I}]$.
(In Figure~\ref{whirlpaths}, $x$ is $1$.)  Let the cyclic permutation
$\Sigma$ of $\mathcal{I}$ be $(I_1,I_2,\ldots,I_{k-1},I_k,\ldots,I_r)$
where the intervals $I_j$ with $x\in I_j$ and $x\ne f_{I_j}$ are
$I_1,I_2,\ldots,I_{k-1}$.  Note that the linear order
$I_1,I_2,\ldots,I_r$, which plays an important role below, has been
specified unless $k$ is $1$ or $k-1$ is $r$.  For $k=1$, let $I_1$ be
the interval $I$ in $\mathcal{I}$ that minimizes the size of the
interval $[x,f_I]$.  For $k-1=r$, let $I_1$ be the interval $I$ in
$\mathcal{I}$ that minimizes the size of the interval $[x,l_I]$.
Consider the subsets $\{p_1,p_2,\ldots,p_k\}$ and
$\{p'_1,p'_2,\ldots,p'_k\}$ of $\mathbb{Z}^2$ where $p_i=(k-i,i-1)$
and $p'_i=p_i+(m,r)$.  Let $L$ and $L'$ be the lines of slope $-1$
that contain these sets.  Let $P$ be the lattice path from $p_1$ to
$p'_1$ formed from the sequence
$x,\sigma(x),\sigma^2(x),\ldots,\sigma^{-1}(x)$ by replacing each
element $l_{I_j}$, for $I_j\in\mathcal{I}$, by a North step and
replacing the other elements by East steps.  Let $Q$ be the lattice
path from $p_k$ to $p'_k$ formed from
$x,\sigma(x),\sigma^2(x),\ldots,\sigma^{-1}(x)$ by replacing each
element $f_{I_j}$, for $I_j\in\mathcal{I}$, by a North step and
replacing the other elements by East steps.  Note that $P$ never goes
above $Q$.  The lines $L$ and $L'$ and the paths $P$ and $Q$ bound the
region of interest.  Label the North and East steps in this region as
follows: steps that are adjacent to the points $p_1,p_2,\ldots,p_k$
are labelled $x$, those one step away from $p_1,p_2,\ldots,p_k$ are
labelled $\sigma(x)$, and so on.  The resulting diagram, which we
denote by $D(\mathcal{I},x)$, depends on both $\mathcal{I}$ and $x$.
(To simplify the example, the diagram shown in Figure~\ref{whirlpaths}
omits the labels on the East steps.)  The diagram $D(\mathcal{I},x)$
captures the set system $\mathcal{I}$: each interval among
$I_k,I_{k+1},\ldots,I_r$ is the set of labels on the North steps in
one row; each interval $I_i$ among $I_1,I_2,\ldots,I_{k-1}$ also
appears in this way, but split into two parts, with $x$ and the
elements in the last part of $I_i-x$ appearing among the lowest $k-1$
rows and with the elements in the first part of $I_i-x$ appearing
among the highest $k-1$ rows.  Theorem~\ref{thm:altrep}, which is a
counterpart of Theorem~\ref{thm:bases}, shows the significance of
$D(\mathcal{I},x)$.

\begin{thm}\label{thm:altrep}
  Fix an element $x$ in a multi-path matroid $M[\mathcal{I}]$. A set
  $B$ is a basis of $M[\mathcal{I}]$ if and only if there is a lattice
  path $R$ such that
\begin{itemize}
\item[(i)] $R$ goes from a point $p_i$ to the corresponding point
  $p'_i$,
\item[(ii)] $R$ uses East and North steps of the diagram
  $D(\mathcal{I},x)$, and
\item[(iii)] the labels on the North steps of $R$ are the elements of
  $B$.
\end{itemize}
\end{thm}

\begin{proof}
  Let $b_1,b_2,\ldots,b_r$, in this order, be the labels on the North
  steps of a path $R$ that satisfies conditions (i) and (ii). Thus,
  $b_1,b_2,\ldots,b_r$ are contained, respectively, in $r$ consecutive
  intervals in the cycle $(I_1,I_2,\ldots,I_r)$; also,
  $b_1,b_2,\ldots,b_r$ are distinct since
  the North and East steps of $R$ are labelled, in order,
  $x,\sigma(x),\sigma^2(x),\ldots,\sigma^{-1}(x)$. It follows that
  $\{b_1,b_2,\ldots,b_r\}$ is a basis of $M[\mathcal{I}]$.

  For the converse, we use the notation established when defining the
  diagram $D(\mathcal{I},x)$.  All references to an order on the
  ground set $S$ are to the linear order
  $x<\sigma(x)<\sigma^2(x)<\cdots<\sigma^{-1}(x)$.  Assume $B$ is a
  basis of $M[\mathcal{I}]$ and let $\phi:B\longrightarrow\mathcal{I}$
  be a matching.  To complete the proof, it suffices prove the
  following claim.
  \begin{quote}
    \emph{The elements of $B$, listed in order as
      $b_1,b_2,\ldots,b_r$, are in the sets
      $I_{k-t},I_{k-t+1},\ldots,I_{k-1},I_k,
      \ldots,I_r,I_1,I_2,\ldots,I_{k-t-1}$, respectively, for some $t$
      with $0\leq t\leq k-1$.}
  \end{quote}
  Indeed, the required path $R$ takes East steps from $p_{k-t}$ until
  a North step labelled $b_1$ is reached; after taking that North
  step, East steps are taken until a North step labelled $b_2$ is
  reached, and so on.
  
  We prove the claim by first constructing a matching for a different
  set system.  For $i$ with $1\leq i \leq k-1$, let $X_i$ be the first
  part of $I_i$ with respect to $x$ and let $Y_i$ be $I_i-X_i$.  Note
  that the set $\phi^{-1}(\{I_1,I_2,\ldots,I_{k-1}\})$ is the disjoint
  union of two subsets whose elements are, in order, say,
  $b'_1,b'_2,\ldots,b'_t$ and $b''_1,b''_2,\ldots,b''_{k-1-t}$, where
  each $b'_i$ is in the subset $Y_j$ of the set $I_j=\phi(b'_i)$ while
  each $b''_i$ is in the subset $X_j$ of the set $I_j=\phi(b''_i)$.
  Thus, $0\leq t\leq k-1$.  Let $\mathcal{I}'$ be the set system that
  consists of the intervals
  $$Y_{k-t},Y_{k-t+1},\ldots,Y_{k-1},I_k,
  \ldots,I_r,X_1,X_2,\ldots,X_{k-t-1}.$$
  We also let
  $Z_1,Z_2,\ldots,Z_r$, respectively, denote these intervals.  Let
  $\Phi:B\longrightarrow \mathcal{I'}$ be given by
  \begin{eqnarray}
    & &\Phi(b) =
   \left\{
     \begin{array}{ll}
      Y_{k-1-t+i}, &\mbox{if $b=b'_i$ with $1\leq i \leq t$;}\\
      X_i, &\mbox{if $b=b''_i$ with $1\leq i \leq k-1-t$;}\\
      \phi(b), &\mbox{if
          $b\in\phi^{-1}(\{I_k,I_{k+1},\ldots,I_r\})$.}
    \end{array}
   \right.  \notag
  \end{eqnarray}
  The inclusions $Y_1\subset Y_2\subset \cdots \subset Y_{k-1}$ and
  $X_{k-1}\subset X_{k-2}\subset \cdots \subset X_1$ imply that $\Phi$
  is a matching.
  
  Finally, to prove the claim it suffices to show that the $i$-th
  element $b_i$ of $B$ is in $Z_i$.  If this statement were false,
  then either $b_i<f_{Z_i}$ or $b_i>l_{Z_i}$.  The first option would
  imply that the $i$ elements $b_1,b_2,\ldots,b_i$ can be in only
  $i-1$ sets, namely $Z_1,Z_2,\ldots,Z_{i-1}$; the second option would
  imply that the $r-i+1$ elements $b_i,b_{i+1},\ldots,b_r$ can be in
  only $r-i$ sets, namely $Z_{i+1},Z_{i+2},\ldots,Z_r$.  Both
  conclusions are contradicted by the matching $\Phi$, so the claim
  and the theorem follow.
\end{proof}

Unlike Theorem~\ref{thm:bases}, the correspondence between paths and
bases in Theorem~\ref{thm:altrep} is not bijective.  For example, the
two paths in Figure~\ref{whirlpaths} indicated by thick lines
correspond to the basis $\{1,3,5\}$. While some bases correspond to a
single path, in general each basis corresponds to a family of paths
that arise from a single word in the alphabet $\{E,N\}$ but starting
at different points among $p_1,p_2,\ldots,p_k$.

\begin{figure}
\begin{center}
\includegraphics[width = 3.75truein]{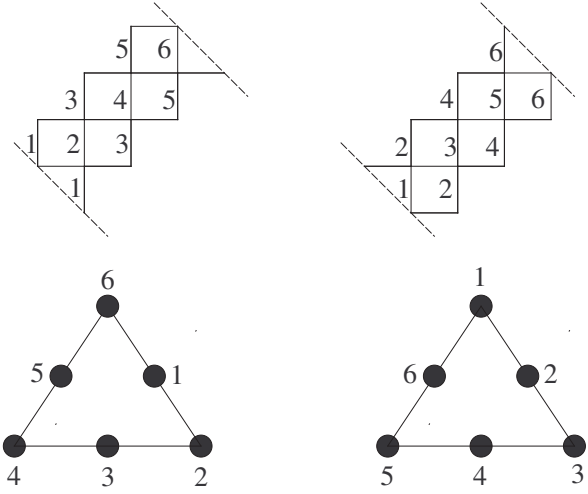}
\end{center}
\caption{The dual of the $3$-whirl $\mathcal{W}^3$ via flipping
diagrams about the line $y=x$.}\label{dual}
\end{figure}

Note that rotating the diagram $D(\mathcal{I},x)$ by $180^\circ$ about
the point $(\frac{m+k-1}{2},\frac{r+k-1}{2})$ gives the diagram
$D\bigl(\mathcal{I},\sigma^{-1}(x)\bigr)$, using the cycle
$\sigma^{-1}$ in place of $\sigma$.

\begin{figure}
\begin{center}
\includegraphics[width = 4.25truein]{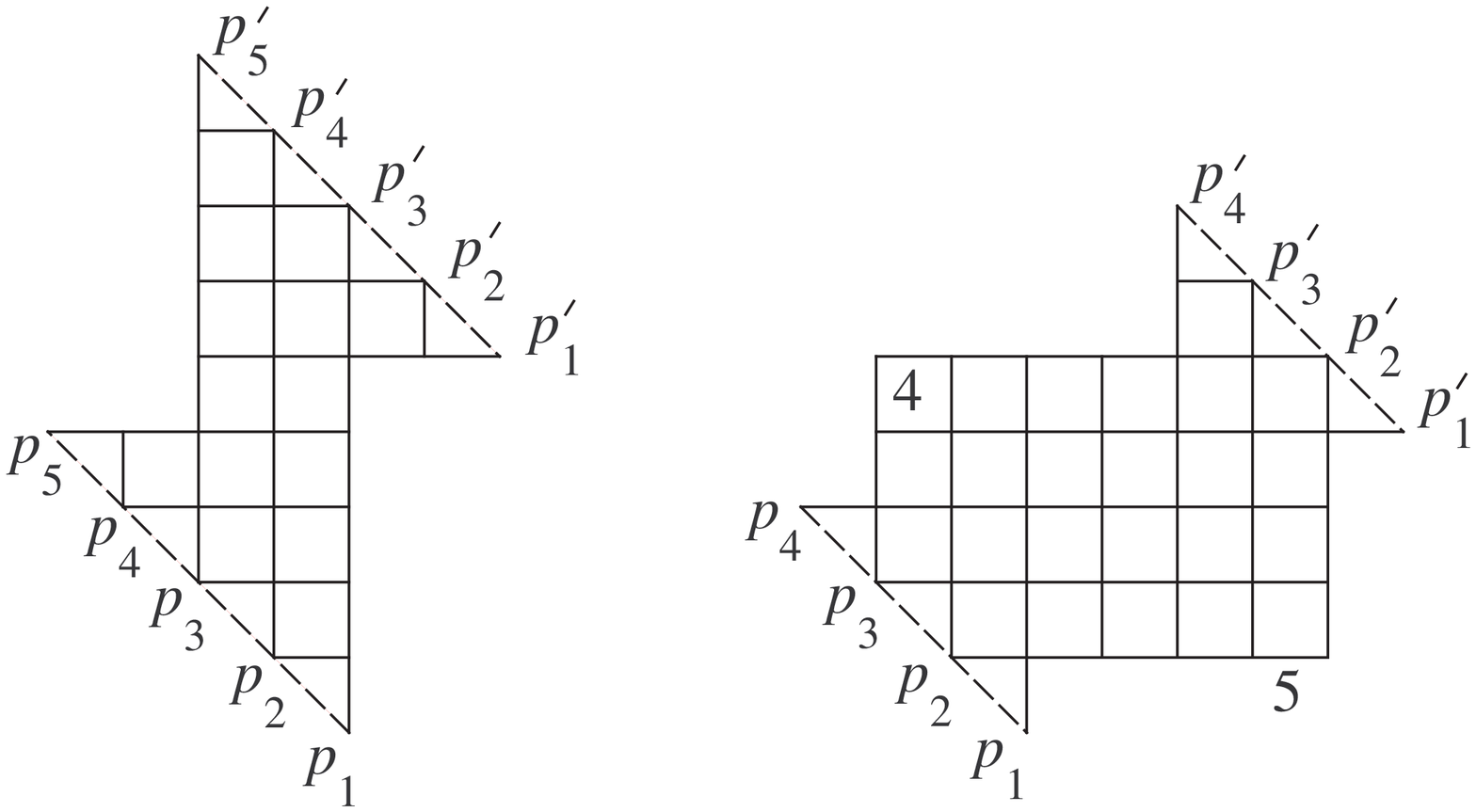}
\end{center}
\caption{Two diagrams $D(\mathcal{I},x)$ that, after reflection in the
  line $y=x$, give set systems that are not
  antichains.}\label{problems}
\end{figure}

Reflecting the diagram $D(\mathcal{I},x)$ in the line $y=x$
interchanges the East and North steps.  Let $D^*(\mathcal{I},x)$
denote this reflected diagram.  A set $X$ is the set of labels on the
the North steps of a path in $D^*(\mathcal{I},x)$ that satisfies
conditions (i) and (ii) of Theorem~\ref{thm:altrep} if and only if $X$
is the complement of a basis of $M[\mathcal{I}]$.  Thus, as
illustrated in Figure~\ref{dual}, $D^*(\mathcal{I},x)$ is a lattice
path representation of the dual matroid $M^*[\mathcal{I}]$.  Some
argument is required, however, to show that $M^*[\mathcal{I}]$ is a
multi-path matroid since the set of $\sigma$-intervals one obtains
from $D^*(\mathcal{I},x)$ need not be an antichain; in particular, the
ground set $S$ may be among these $\sigma$-intervals.  For instance,
the East steps of a column of $D(\mathcal{I},x)$ (for example, the
column between $p_1$ and $p_2$, or that between $p_2$ and $p_3$ in the
first diagram in Figure~\ref{problems}) may be labelled with all
elements of $S$.  Also, the first part of an interval that includes
$x$, say between $p_i$ and $p_{i+1}$, must be joined with with the
corresponding last part between $p'_i$ and $p'_{i+1}$, and this union
may be $S$; the second diagram in Figure~\ref{problems} illustrates
this point with the column between $p_2$ and $p_3$ (the last of the
East steps, labelled $1,2,3,4$, is marked) and that between $p'_2$ and
$p'_3$ (the first of the East steps, labelled $5,6,7,8,9$, is marked).
One way to address this problem, in the spirit of the proofs of
Theorems~\ref{thm:del} and~\ref{thm:con}, is to show how to modify the
set system that corresponds to $D^*(\mathcal{I},x)$ to obtain a
presentation of $M^*[\mathcal{I}]$ by an antichain of
$\sigma$-intervals.  Instead, we introduce a more general type of
diagram (which plays a key role in Section~\ref{sec:tutte}) and show
that for such a diagram $D$, the sets of labels of the North steps of
the paths of $D$ that satisfy conditions (i) and (ii) of
Theorem~\ref{thm:altrep} are the bases of a multi-path matroid.  To
avoid excess terminology, we also call these more general objects,
which we define below, diagrams; this should create no confusion.

A \emph{diagram} $D$ is a $5$-tuple
$(k,m,r,P,Q)$, where $k$ is a positive integer, $m$ and $r$ are
non-negative integers, $P$ is a lattice path from $(k-1,0)$ to
$(k-1+m,r)$, and $Q$ is a lattice path from $(0,k-1)$ to $(m,k-1+r)$
that never goes below $P$.  For $i$ with $1\leq i\leq k$, let $p_i$ be
$(k-i,i-1)$ and let $p'_i$ be $p_i+(m,r)$.  Let $L$ and $L'$ be the
lines of slope $-1$ that contain the points $p_i$ and $p'_i$,
respectively.  The \emph{region} $R(D)$ of a diagram $D$ is the set of
points in $\mathbb{R}^2$, including the boundary, enclosed by the
paths $P$ and $Q$ and the lines $L$ and $L'$.  The \emph{edges} of $D$
are the segments between lattice points in $D$ that are distance $1$
apart.  Assign label $i$ to an edge in $D$ if it is the $i$-th step in
some lattice path that starts at a point on $L$; thus, edges are
labelled with the elements of $[m+r]$.  A \emph{b-path} is a lattice
path contained in the region $R(D)$ that starts at a point $p_i$ and
ends in the corresponding point $p_i'$. Thus, any b-path contains $r$
North steps and $m$ East steps, and the edges are labelled, in order,
$1,2,\ldots,m+r$.  The \emph{label-set} of a b-path $T$ is the set of
labels on the North steps of $T$. Let $\mathcal{B}(D)$ be the set of
all label-sets of b-paths in $D$.  We now show that $\mathcal{B}(D)$
is the set of bases of a multi-path matroid, which we denote by
$M[D]$.  (To recover multi-path matroids in complete generality,
replace the labels $1,2,\ldots,m+r$ with the elements
$x,\sigma(x),\ldots,\sigma^{-1}(x)$, respectively.  In much of the
rest of the paper, we favor the notational simplicity gained by having
$[m+r]$ be the ground set of $M[D]$.)

\begin{thm}\label{thm:gendia}
  For any diagram $D= (k,m,r,P,Q)$, the collection $\mathcal{B}(D)$ of
  subsets of $[m+r]$ is the set of bases of a multi-path matroid.
\end{thm}

\begin{figure}
\begin{center}
\includegraphics[width = 3.0truein]{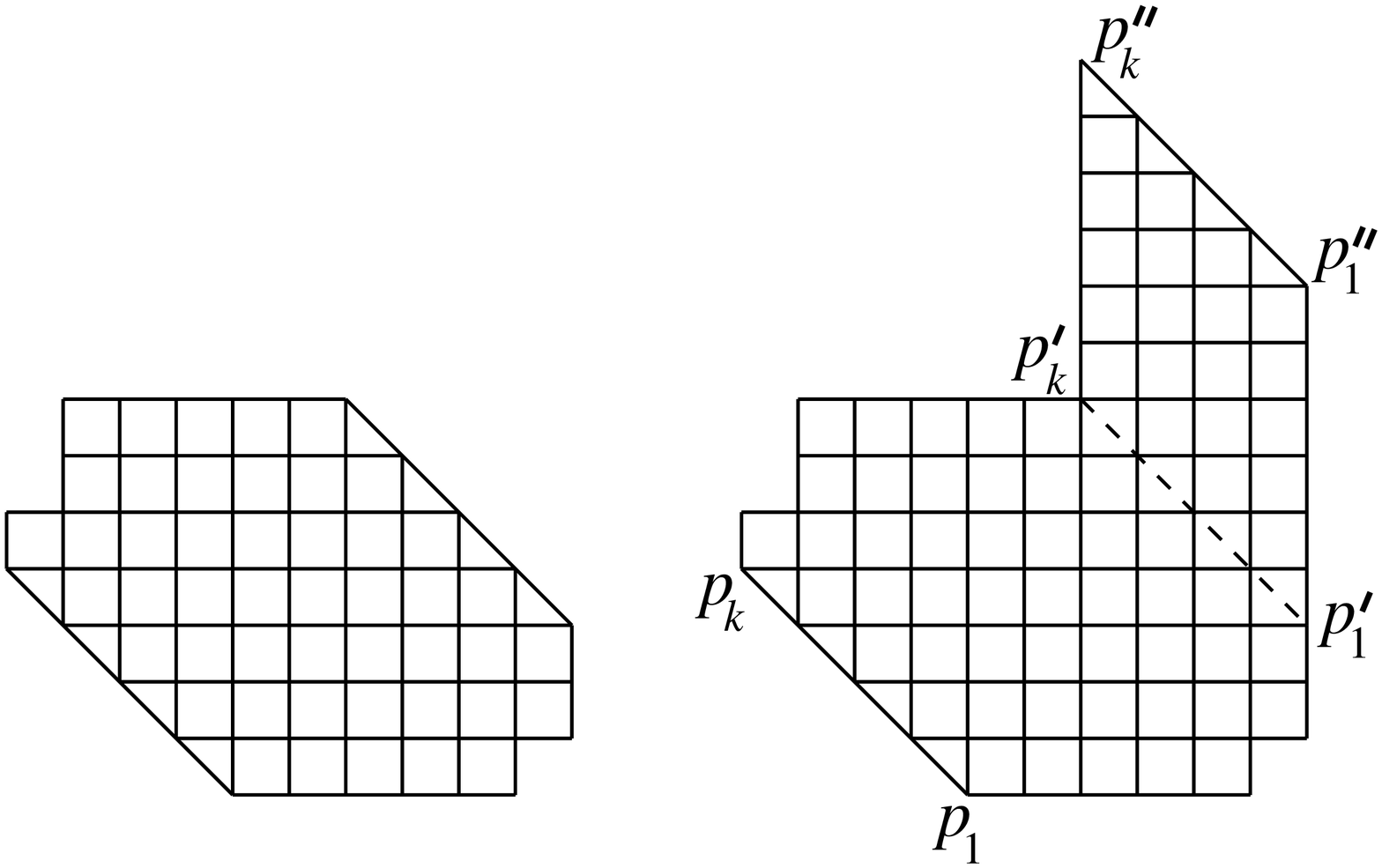}
\end{center}
\caption{The diagram $D=(5,6,3,E^5NEN^2,NEN^2E^5)$ and its
extension $D'=(5,6,9,E^5NEN^8,NEN^2E^5N^6)$.)}\label{extend}
\end{figure}

\begin{proof}
  By Theorem~\ref{thm:con}, it suffices to prove that $\mathcal{B}(D)$
  is the set of bases of a contraction of a multi-path matroid
  $M[\mathcal{I}]$.  Toward this end, let $D'$ be the diagram
  $(k,m,r+k+1,PN^{k+1},QN^{k+1})$.  (See Figure~\ref{extend}.)  For
  $i$ with $1\leq i\leq k$, let $p_i$ and $p'_i$ be as above and let
  $p''_i$ be $p'_i+(0,k+1)$.  Denote the rows of $D'$, from the bottom
  up, by $R_1,R_2,\ldots,R_{r+2k}$.  Let $\sigma$ be the cycle
  $(1,2,\ldots,m+r+k+1)$.  Let $\mathcal{I}$ consist of the following
  sets: $I_j$, for $j$ with $1\leq j<k$, is the union of the set of
  labels on the North steps of row $R_j$ and that of row
  $R_{r+k+j+1}$; the set $I_j$, for $j$ with $k\leq j\leq k+r+1$,
  consists of the labels on the North steps in row $R_j$.  Each set
  $I_j$ is a $\sigma$-interval and $D'$ is the diagram
  $D(\mathcal{I},1)$ for $\mathcal{I}$.  We claim that $\mathcal{I}$
  is an antichain.  Note that each set $I_j$ has at most $m+k$
  elements and so is a proper subset of $[m+r+k+1]$.  The sets
  $I_k,I_{k+1},\ldots,I_{k+r+1}$ form an antichain since we have
  $f_{I_k}<f_{I_{k+1}}<\cdots<f_{I_{k+r+1}}$ and
  $l_{I_k}<l_{I_{k+1}}<\cdots<l_{I_{k+r+1}}$ for these intervals in
  the usual linear order on $[m+r+k+1]$.  A similar argument shows
  that $I_1,I_2,\ldots,I_{k-1}$ form an antichain. Now consider $I_h$
  and $I_j$ with $1\leq h<k\leq j\leq r+k+1$.  At least one of $1$ and
  $m+r+k+1$ is not in $I_j$, so $I_h\not \subseteq I_j$. The
  containment $I_j\subset I_h$ would imply that either $I_j\subseteq
  [f_{I_h},m+r+k+1]$ or $I_j\subseteq [1,l_{I_h}]$ holds.  The first
  inclusion contradicts the inequality $1\leq f_{I_j}<f_{I_h}\leq
  m+r+k+1$ that is evident from the diagram $D'$; the second
  containment contradicts the inequality $1\leq l_{I_h}<l_{I_j}\leq
  m+r+k+1$ that is also evident from $D'$.  Thus, $\mathcal{I}$ is an
  antichain of $\sigma$-intervals.

  Let $Z$ consist of the last $k+1$ elements of $[m+r+k+1]$.  We now
  show that $\mathcal{B}(D)$ is the set of bases of the contraction of
  the multi-path matroid $M[\mathcal{I}]$ by $Z$.  Since $Z$ is
  independent in $M[\mathcal{I}]$, the bases of $M[\mathcal{I}]/Z$ are
  the subsets $B$ of $[m+r]$ for which $B\cup Z$ is a basis of
  $M[\mathcal{I}]$.  Note that the last $k+1$ steps in any lattice
  path whose label set is $B\cup Z$ are North steps that go from a
  point $p'_i$ to the corresponding point $p''_i$.  Thus, $B\cup Z$ is
  a basis of $M[\mathcal{I}]$ if and only if $B\cup Z$ is the label
  set of a path in $D'$ that goes from some point $p_i$ to the
  corresponding point $p''_i$ through the point $p'_i$.  It follows
  that $B$ is a basis of $M[\mathcal{I}]/Z$ if and only if $B$ is the
  label set of a b-path in $D$, that is, if and only if $B$ is in
  $\mathcal{B}(D)$, as claimed.
\end{proof}

That arbitrary diagrams define multi-path matroids allows us to give
another perspective on certain minors.  (Since the proof of
Theorem~\ref{thm:gendia} uses Theorem~\ref{thm:con}, this does not
replace our earlier work.)  Let $M$ be the multi-path matroid on
$[m+r]$ that is represented by a diagram $D$.  Let $X$ and $Y$ be
disjoint subsets of $[m+r]$ where $Y$ is independent, $X$ is
coindependent (i.e., the complement of a spanning set), and $X\cup Y$
consists of the last $k$ elements of $[m+r]$.  From the formulation of
minors in terms of bases, it follows that the bases of the minor
$M\del X/Y$ correspond to the paths in $D$ whose last $k$ steps are
determined: these steps are East or North according to whether their
labels are in $X$ or $Y$, respectively.  The initial segments of paths
in $D$ whose last $k$ steps are as specified make up a smaller diagram
$D'$ that represents the minor $M\del X/Y$.  This observation, which
is behind the proof of Theorem~\ref{thm:gendia}, plays an important
role in the next section.

The next theorem summarizes the results in this section. The assertion
about duality follows from Theorem~\ref{thm:gendia} and the remarks
before that theorem.

\begin{thm}\label{thm:minordual}
  The class of multi-path matroids is dual-closed, minor-closed, and
  properly contains the class of lattice path matroids.
\end{thm}

\section{Tutte Polynomial}\label{sec:tutte}

The Tutte polynomial has received considerable attention, in part due
to its many striking properties (e.g., it is the universal
deletion-contraction invariant) and its many important evaluations
(e.g., the chromatic and flow polynomials of a graph, the weight
enumerator of a linear code, and the Jones polynomial of an
alternating knot). (See~\cite{bo,ckcc}.)  In this section, we show
that the Tutte polynomial of a multi-path matroid can be computed in
polynomial time.  This result stands in contrast to the hardness
results known for computing the Tutte polynomial of an arbitrary
member of many classes of matroids~\cite{cpv,bc,hard,bic,bip}.  We
cast our work on the Tutte polynomial in a broader framework; we
introduce what we call computation graphs, which allow us to apply
dynamic programming.

The Tutte polynomial $t(M;x,y)$ of a matroid $M$ on the ground set
$S$ can be defined in a variety of ways, perhaps the most basic
of which is the following:
\begin{equation}\label{tdef}
t(M;x,y) = \sum_{A \subseteq S} (x-1)^{r(S)-r(A)}(y-1)^{|A|-r(A)}.
\end{equation}
The following recurrence relation is more suited to our work.  The
Tutte polynomial $t(M;x,y)$ is $1$ if $M$ is the empty matroid;
otherwise, for any element $e$ of $M$,
\begin{eqnarray} \label{rdef}
 & &t(M;x,y)=
   \left\{
     \begin{array}{ll}
      x\,t(M/e;x,y) &\mbox{if $e$ is an isthmus;}\\
      y\,t(M\backslash e;x,y) &\mbox{if $e$ is a loop;}\\
      t(M/e;x,y)+t(M\backslash e;x,y) &\mbox{otherwise.}
    \end{array}
   \right.
\end{eqnarray}

As stated, both of these formulations require roughly $2^{|S|}$
computations.  We take advantage of the fact that for a multi-path
matroid the recurrence relation can be applied in a manner that
involves minors that are easily recognized to be equal; more
precisely, the number of different minors that need to be considered
turns out to be polynomial in $|S|$, and this allows us to organize
the computation in a way that runs in polynomial time.  Before turning
to multi-path matroids, we establish a general framework for
computations of this type.

Let $M$ be a matroid on the set $[n]$. To use the recurrence
relation~(\ref{rdef}), it suffices to consider what we will call the
\emph{initial minors} of $M$, that is, the matroids formed by deleting
or contracting, in turn, $n, n-1, \ldots, h+2, h+1$, where at the
stage at which an element is deleted, it is not an isthmus, and at the
stage at which an element is contracted, it is not a loop. The ground
set of an initial minor is an initial segment $[h]$ of $[n]$. Note
that if $M\backslash X/Y$ is an initial minor, then $Y$ is independent
and $X$ is coindependent.

We define a \emph{computation graph} $G$ for the matroid $M$ to be an
edge-labelled directed graph with label set $\{c,d\}$ that satisfies
the following conditions.
\begin{itemize}
\item[(1)] Each vertex $u$ represents an initial minor $M_u$ of $M$.
  Every initial minor of $M$ is represented by at least one vertex.
\item[(2)] Let $u$ be a vertex and let $h$ be the greatest element of
  the initial minor $M_u$. If there is a $d$-edge from $u$ to $v$,
  then $M_u\backslash h=M_v$. If there is a $c$-edge from $u$ to $w$,
  then $M_u/h=M_w$. In addition,
\begin{itemize}
\item[(a)] if $h$ is an isthmus of $M_u$, then $u$ is the tail of
  exactly one $c$-edge and no $d$-edge;
\item[(b)] if $h$ is a loop of $M_u$, then $u$ is the tail of exactly
  one $d$-edge and no $c$-edge;
\item[(c)] otherwise $u$ is the tail of exactly one $c$-edge and one
  $d$-edge.
\end{itemize}
\item[(3)] There are two distinguished vertices $v_M$ and
  $v_\emptyset$; these are the unique vertices that represent the
  matroid $M$ and the empty matroid, respectively.
\end{itemize}

By point (1), to construct a computation graph $G$ for a matroid $M$
by using some representation (e.g., a multi-path diagram), apart from
the trivial cases in point (3) we are not required to determine
whether different representations give the same minor.  Note that the
restrictions imposed on the edges imply that $u$ is at distance $h$
from $v_\emptyset$ if and only if $M_u$ has $h$ elements; let $V_h$ be
the set of such vertices $u$. Then $\{V_0, V_1, \ldots, V_n\}$ is a
partition of the vertices of $G$ and any edge that has its tail in
$V_h$ has its head in $V_{h-1}$.

Recurrence relation~(\ref{rdef}) allows us to compute $t(M;x,y)$ from
the computation graph $G$. There is a trade-off between several
factors that enter into the computation graph: having fewer vertices
allows us to compute the Tutte polynomial more quickly, but getting
fewer vertices requires recognizing that many initial minors (perhaps
with different representations) are equal.  A typical application of
these ideas would yield a computation graph with polynomially many
vertices without determining all instances of equal initial minors.
The following lemma helps quantify these observations.

\begin{lemma}\label{lemma:compgraph}
  We can compute the Tutte polynomial $t(M;x,y)$ from a computation
  graph $G$ on $\nu$ vertices in $\mathcal{O}(\nu rm)$ operations,
  where $r$ and $m$ are the rank and nullity of $M$.
\end{lemma}

\begin{proof}
  Partition the vertices of $G$ into blocks $V_0,\ldots, V_{m+r}$, as
  described before; since $G$ has no oriented cycles this can be done
  with $\mathcal{O}(\nu)$ operations. Assign to every vertex $u$ the
  Tutte polynomial $t(M_u;x,y)$ in the following manner. First assign
  $1$ to the unique vertex $v_\emptyset$ in $V_0$, then compute the
  Tutte polynomials for all vertices in $V_1$, then those for all
  vertices of $V_2$, and so on. To compute the Tutte polynomial
  $t(M_u;x,y)$ for $u$ in $V_h$, apply recurrence
  relation~(\ref{rdef}): by condition (2) in the definition of a
  computation graph, the edges for which $u$ is the tail indicate
  which of the three cases of the recurrence to use, and the Tutte
  polynomials of $M_u\backslash h$ and $M_u/h$ have already been
  computed because they correspond to vertices of $V_{h-1}$. Thus for
  every vertex $u$ we just need to add two polynomials or multiply a
  polynomial by $x$ or $y$, and this can be done in $\mathcal{O}(rm)$
  operations since $t(M;x,y)$ has at most $rm+r+m$ coefficients.
  Hence we can compute the Tutte polynomial of every initial minor of
  $M$, including $M$ itself, in $\mathcal{O}(\nu+\nu rm)$, that is,
  $\mathcal{O}(\nu rm)$, operations.
\end{proof}

We now focus on the multi-path matroid $M[\mathcal{I}]$, or $M$, on
$[m+r]$ where $\sigma$ is the cycle $(1,2,\ldots,m+r)$.  Let $D$ be
the diagram $D(\mathcal{I},1)$.  We first study the initial minors
$M\backslash X/Y$ that arise in constructing a computation graph for
$M$, and to do so we work with the diagrams introduced in
Section~\ref{sec:basics}.  In particular, we show how to obtain a
diagram $D'$ for any initial minor $M\backslash X/Y$. The resulting
diagrams need not arise from $\sigma$-interval presentations.

\begin{figure}
\begin{center} 
\includegraphics[width = 3.25truein]{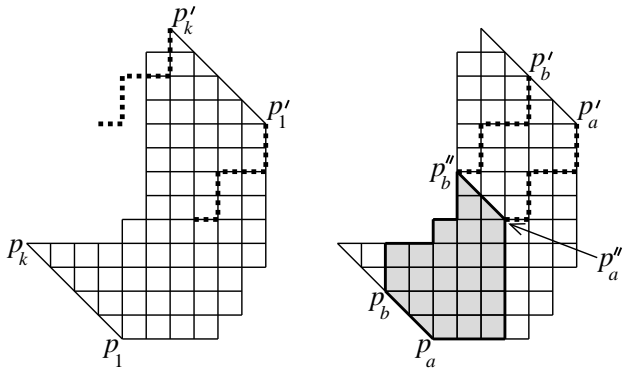}
\end{center}
\caption{The shaded region in the second diagram represents the initial 
minor
  $M\del \{15,14,11,10\}/\{13,12,9\}$.}
\label{fig:initial-minors-diagram}
\end{figure}

It follows from the basis formulation of deletion and contraction that
$B$ is a basis of $M\backslash X/Y$ if and only if $B\cup Y$ is a
basis of $M$ (recall that $X$ and $Y$ are, respectively, coindependent
and independent). These bases, by Theorem~\ref{thm:altrep},
corresponds to b-paths where the last $q=|X\cup Y|$ steps are
determined: steps corresponding to elements of $Y$ are North and steps
corresponding to elements of $X$ are East. Let $a$ and $b$ be the
smallest and largest integers $i$ such that there is a path from $p_i$
to $p_i'$ in $D$ with the last $q$ steps as specified by $X$ and $Y$.
(See Figure~\ref{fig:initial-minors-diagram}.) For $i$ between $a$ and
$b$ let $p_i''$ be the point $p_i'-(|X|,|Y|)$; thus any path from
$p_i$ to $p_i'$ whose last $q$ steps are as specified by $X$ and $Y$
goes through the point $p_i''$.  Let $P'$ be the lattice path in $D$
from $p_a$ to $p_a''$ that no path in $D$ from $p_a$ to $p_a''$ goes
below; similarly, let $Q'$ be the lattice path in $D$ from $p_b$ to
$p_b''$ that no path in $D$ from $p_b$ to $p_b''$ goes above. Let $D'$
be the diagram that has $p_a, \ldots, p_b$ as starting points,
$p_a'',\ldots, p_b''$ as ending points, and $P'$ and $Q'$ as the
bottom and top border. Thus, if $D$ is $(k,m,r,P,Q)$, then $D'$ is
$(b-a+1, m-|X|, r-|Y|, P', Q')$.

\begin{lemma}\label{lemma:minors}
  Let $D=(k,m,r,P,Q)$ be the diagram $D(\mathcal{I}, 1)$ of a
  multi-path matroid $M$ on $[n]$.
\begin{itemize}
\item[(1)] We can construct from $D$ a diagram $D'$ corresponding to
  an initial minor $M\backslash X/Y$ in $\mathcal{O}(n)$ operations.
\item[(2)] We can construct from $D$ at most
  $(n+1)\bigl(\text{min}(r,m)+1\bigr)(k^2+k)/2$ different diagrams
  $D'$ corresponding to initial minors of $M$. In particular, $M$ has
  at most this many initial minors.
\end{itemize}
\end{lemma}

\begin{proof}
  The description above for constructing $D'$ from $D$ has two parts:
  find $a$ and $b$, and then find $P'$ and $Q'$. We sketch how to do
  these two steps.  Since $X$ and $Y$ are coindependent and
  independent, $a$ and $b$ exist; find them by comparing the last
  $|X\cup Y|$ steps of $P$ and $Q$ with the steps specified by $X$ and
  $Y$. (See Figure~\ref{fig:initial-minors-diagram}. The dotted paths
  are those specified by $X$ and $Y$.)  Specifically: let
  $\textrm{N}(W,i)$ be the number of North steps among the last $i$
  steps of a path $W$ and let $P_{X,Y}$ be the path specified by $X$
  and $Y$; then $a$ is $$\max\{ \textrm{N}(P_{X,Y},i)-\textrm{N}(P,i)
  :\, 0\leq i \leq |X\cup Y| \}+1.$$
  A similar formula gives $b$, so
  $a$ and $b$ can be computed in $\mathcal{O}(n)$ operations.
  Construct $P'$ (respectively, $Q'$) by going from $p_a$ to $p_a''$
  (respectively, from $p_b$ to $p_b''$), taking East (respectively,
  North) steps whenever possible. This also takes $\mathcal{O}(n)$
  operations.

  Assertion (2) follows by noting that a diagram is completely
  determined by (i) the size of $X\cup Y$, (ii) the size of either $X$
  and $Y$, (iii) the points $p_a''$ and $p_b''$, and that these two
  points are determined by the two numbers $a \leq b$ between $1$ and
  $k$.
\end{proof}

We now show how to compute the Tutte polynomial of a multi-path
matroid $M$ in polynomial time from its diagram $D=D(\mathcal{I},1)$.
We start by constructing a computation graph for $M$ whose vertices
correspond to the diagrams of initial minors of $M$ that are obtained
from $D$ as described before Lemma~\ref{lemma:minors}. Start with a
graph that consists of just one vertex $v_M$ that corresponds to the
diagram $D$ and iterate the following process.
\begin{itemize}
\item Choose a vertex $v$ other than $v_\emptyset$ with outdegree 0.
  Let $M'$, on $[h]$, and $D'$ be the corresponding initial minor and
  diagram.
\item Compute the diagrams $D_d$ and $D_c$ corresponding to
  $M'\backslash h$, if $h$ is not an isthmus, and $M'/h$, if $h$ is
  not a loop, as described before Lemma~\ref{lemma:minors}.
\item Find a vertex $v_d$ (respectively, $v_c$) in the computation
  graph that corresponds to $D_d$ (respectively, $D_c$); if there is
  no such vertex, add a new vertex to the computation graph. Add a
  $d$-edge (respectively, a $c$-edge) from $v$ to $v_d$ (respectively,
  $v_c$).
\end{itemize}
Stop when the only vertex of outdegree 0 is $v_\emptyset$, which
corresponds to the empty matroid. The resulting graph $G$ is clearly a
computation graph for $M$.  The same initial minor $M'$ can be
represented more than once in $G$ since different diagrams can
represent it, but each diagram $D'$ appears just once and all diagrams
have been derived from $D$.  By part (2) of Lemma~\ref{lemma:minors},
the number $\nu$ of vertices of $G$ is
$\mathcal{O}(n\,\text{min}(r,m)k^2)$.  By Lemma~\ref{lemma:compgraph},
we can compute $t(M;x,y)$ from $G$ in $\mathcal{O}(rm\nu)$ operations.
So now we need only show that this construction of the computation
graph can be done in polynomial time.

We show that we can construct $G$ in $\mathcal{O}(\nu\, n\log \nu)$
operations. Consider the operations required for each iteration of the
algorithm (each expansion of a vertex $v$ of outdegree 0). First we
compute $D_d$ and $D_c$ in $\mathcal{O}(n)$ operations and then we
check whether they are already in the graph.  Comparing two diagrams
(i.e., $5$-tuples) requires $\mathcal{O}(n)$ operations; by using a
suitable ordering of the vertices, a binary search using
$\mathcal{O}(\log \nu)$ comparisons suffices to determine whether a
given diagram is already in the graph.  Thus we need
$\mathcal{O}(n\log \nu)$ operations for any of the $\nu$ iterations,
so $G$ can be constructed in $\mathcal{O}(\nu\, n\log \nu)$
operations.

Hence the number of operations needed to construct this
computation graph and obtain the Tutte polynomial from it is
$\mathcal{O}\bigr(\nu(rm+n\log \nu)\bigr)$. To simplify the
expression for the number of operations required, note that $r+m$
is $n$ and $k$ is less than $n$; also, $\log \nu$ is
$\mathcal{O}(\log n)$ because $\nu$ depends polynomially on $n$.
Thus, the work in this section gives the following theorem.

\begin{thm}\label{thm:Tutte-poly}
  We can compute the Tutte polynomial of a multi-path matroid on $n$
  elements in $\mathcal{O}(n^6)$ operations.
\end{thm}

\section{Basis Activities}\label{sec:activities}

Another formulation of the Tutte polynomial is given by basis
activities, which are also of independent interest. In this section,
we describe the internal and external activities of bases of
multi-path matroids in terms of lattice paths in diagrams and we
sketch an alternative approach to computing the Tutte polynomial of a
multi-path matroid through basis activities.

The Tutte polynomial of $M$ can be written as
\begin{equation}\label{tactive} t(M;x,y)= \sum_{B\in
\mathcal{B}(M)} x^{i(B)}y^{e(B)},
\end{equation}
where $\mathcal{B}(M)$ is the collection of bases of $M$ and the
exponents $i(B)$ and $e(B)$ are defined as follows.  Fix a linear
order $<$ on the ground set $S$ of $M$ and let $B$ be a basis of
$M$. An element $u$ in $S-B$ is \emph{externally active with
respect to $B$} if there is no element $v$ in $B$ with $v<u$ for
which $(B-v)\cup u$ is a basis.  An element $u$ in $B$ is
\emph{internally active with respect to $B$} if there is no
element $v$ in $S-B$ with $v<u$ for which $(B-u)\cup v$ is a
basis. The \emph{internal activity} $i(B)$ of a basis $B$ is the
number of elements that are internally active with respect to $B$.
The external activity of $B$, denoted $e(B)$, is defined
similarly. Note that $i(B)$ and $e(B)$ depend not only on $B$ but
also on the order. Equation~(\ref{tactive}) says that the
coefficient of $x^iy^e$ in $t(M;x,y)$ is the number of bases of
$M$ with internal activity $i$ and external activity $e$. In
particular, the number of such bases is independent of the order.

We will use the following lemma, which is well-known and easy to prove.

\begin{lemma}\label{dualact}
  Fix a linear order on the ground set $S$ of a matroid $M$ and its
  dual $M^*$. An element $u$ is internally active with respect to the
  basis $B$ of $M$ if and only if $u$ is externally active with
  respect to the basis $S-B$ of $M^*$.
\end{lemma}

Throughout this section we use the notation and terminology we
establish in the next several paragraphs. We assume that the ground
set of the multi-path matroid $M[\ICal]$ is $[m+r]$ and that $\sigma$
is the cycle $(1,2,\ldots,m+r)$.  We study the internal and external
activities of the bases of $M[\ICal]$ relative to the linear order
$1<2<\cdots<m+r$. Let $D$ be the diagram $D(\ICal, 1)=(k, m, r, P,
Q)$; recall that $P$ and $Q$ are respectively the bottom and top
border of the diagram.

For any subset $X$ of $[m+r]$ the \emph{representation $\Pi(X,p)$ of
  $X$ starting at the lattice point $p$} is the path of $m+r$ steps
that starts at $p$ whose $u$-th step is $\NStep$ if $u$ is in $X$, and
$\EStep$ otherwise.  We say that a path is \emph{valid} if it is
entirely contained in the diagram $D$.  Thus, Theorem~\ref{thm:altrep}
states that the bases of $M[\ICal]$ are the sets $B$ such that, for
some $p_i$, the path $\Pi(B,p_i)$ is valid and ends at the
corresponding point $p_i'$. Note that if $\Pi(B,p_i)$ and $\Pi(B,p_j)$
are both valid paths for $i<j$, then all paths $\Pi(B,p_t)$ with
$i<t<j$ are also valid.

For $v,u$ in $[m+r]$ with $v\leq u+1$, we use $[v,u]\Pi(X,p_i)$ to
denote the path that starts at the beginning of the $v$-th step of
$\Pi(X,p_i)$ and follows this path until the end of the $u$-th step.
The notation $(v,u]\Pi(X,p_i)$, $[v,u)\Pi(X,p_i)$, and
$(v,u)\Pi(X,p_i)$ is defined in the obvious way; for instance
$(v,u)\Pi(X,p_i)$ is $[v+1,u-1]\Pi(X,p_i)$.  In particular,
$(v,u)\Pi(X,p_i)$ is defined when $v\leq u-1$, and 
$(u-1,u)\Pi(X,p_i)$ consists of the single point that is common to
steps $u-1$ and $u$ of $\Pi(X,p_i)$.

The following lemma gives the conditions under which an element $u$ of
a basis $B$ can be replaced by an element $v$ to yield a basis $B'$.

\begin{lemma}\label{lem:basisexchange}
  Let $B$ be a basis of a multi-path matroid $M[\ICal]$ with $u\in B$,
  $v\notin B$. Let $\Pi(B,p_i)$ be a valid path. Let
  $B'$ be $(B-u)\cup v$.
\begin{enumerate}
\item If $v<u$, then $B'$ is a basis if and only if either
\begin{enumerate}
\item the path $(v,u)\Pi(B,p_i)$ does not touch the top border $Q$, or
\item neither $[1,v)\Pi(B,p_i)$ nor $(u,m+r]\Pi(B,p_i)$ touches the 
bottom border $P$.
\end{enumerate}
\item If $u<v$, then $B'$ is a basis if and only if either
\begin{enumerate}
\item the path $(u,v)\Pi(B,p_i)$ does not touch $P$, or
\item neither $[1,u)\Pi(B,p_i)$ nor $(v,m+r]\Pi(B,p_i)$ touches $Q$.
\end{enumerate}
\end{enumerate}
\end{lemma}

\begin{figure}
\begin{center}
\includegraphics[width = 3.0truein]{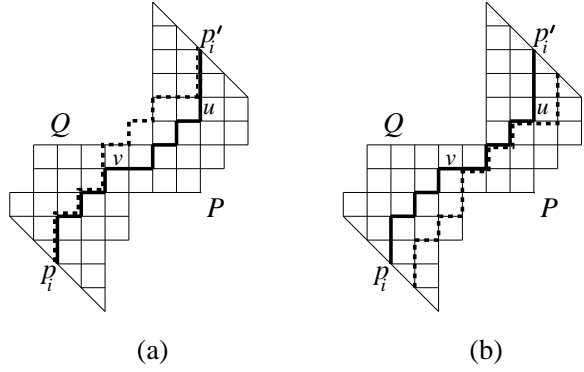}
\end{center}
\caption{Paths $\Pi(B',p_i)$ (dotted line) and $\Pi(B,p_i)$ in
part (a), and $\Pi(B',p_{i-1})$ (dotted line) and $\Pi(B,p_i)$ in
part (b).} \label{fig:basisexchangelemma}
\end{figure}

\begin{proof}
  By duality it suffices to prove the first claim.  By
  Theorem~\ref{thm:altrep}, $B'$ is a basis if and only if $B'$ has a
  valid representation. Compare the paths $\Pi(B',p_i)$ and
  $\Pi(B',p_{i-1})$ with $\Pi(B,p_i)$. (See
  Figure~\ref{fig:basisexchangelemma}.)  Since $\Pi(B',p_i)$ is above
  $\Pi(B,p_i)$, only $Q$ may prevent $\Pi(B',p_i)$ from being valid,
  and in that case no $\Pi(B',p_j)$ with $j\geq i$ would be valid;
  similarly, only $P$ may prevent $\Pi(B',p_{i-1})$ from being valid,
  and in that case no $\Pi(B',p_{j})$ with $j\leq i-1$ would be valid.
  Thus $B'$ is a basis if and only if either $\Pi(B',p_i)$ or
  $\Pi(B',p_{i-1})$ is valid. These two conditions are equivalent to
  conditions (1.a) and (1.b), as Figure~\ref{fig:basisexchangelemma}
  illustrates.
\end{proof}

With the help of this basis exchange lemma we now characterize the
internally and externally active elements.

\begin{thm}\label{thm:activities}
  Let $B$ be a basis of a multi-path matroid $M[\ICal]$ and let
  $\Pi(B,p_i)$ be a valid path.
\begin{enumerate}
\item[(I)] An element $u$ in $B$ is internally active if and only if
  either
\begin{enumerate}
\item $[u] \subseteq B$, or
\item the $u$-th step of $\Pi(B,p_i)$ lies in the top border $Q$ and
  $(u,m+r]\Pi(B,p_i)$ touches the bottom border $P$.
\end{enumerate}
\item[(II)] An element $u$ not in $B$ is externally active if and only 
if
  either
\begin{enumerate}
\item $[u] \cap B = \emptyset$, or
\item the $u$-th step of $\Pi(B,p_i)$ lies in $P$ and
  $(u,m+r]\Pi(B,p_i)$ touches $Q$.
\end{enumerate}
\end{enumerate}
\end{thm}

\begin{proof}
  By duality, we only need to prove part (I).
  
  Note that $u$ is internally active if $[u]\subseteq B$. Thus, let
  $V$ be $[u]-B$ and assume that $V$ is not empty.  Sufficiency
  follows because if $u$ satisfies condition (I.b), then it satisfies
  neither conditions (1.a) nor (1.b) of Lemma~\ref{lem:basisexchange}
  for any $v$ in $V$.
  
  To prove the converse assume that $u$ is internally active. Let $v$
  be $\textrm{max}(V)$. Since $u$ is internally active, $(B-u)\cup v$
  is not a basis, so by condition (1.a) of
  Lemma~\ref{lem:basisexchange} the path $(v,u)\Pi(B,p_i)$ touches
  $Q$.  This path has only North steps, by the choice of $v$, so its
  ending point has to touch $Q$. This proves that the $u$-th step of
  $\Pi(B,p_i)$ lies in $Q$, so the first part of condition (I.b)
  holds.
  
  For the second part, let $v$ be $\textrm{min}(V)$.  By condition
  (1.b) of Lemma~\ref{lem:basisexchange}, since $(B-u)\cup v$ is not a
  basis, at least one of the paths $[1,v)\Pi(B,p_i)$ and
  $(u,m+r]\Pi(B,p_i)$ touches $P$. We show that if the first path
  touches $P$, then so does the second, hence either way
  $(u,m+r]\Pi(B,p_i)$ touches $P$, which proves the second part of
  condition (I.b).  Indeed, the minimality of $v$ implies that
  $[1,v)\Pi(B,p_i)$ has only North steps. So if $[1,v)\Pi(B,p_i)$
  touches $P$, then it has to touch it from the beginning, that is,
  $p_i$ has to be $p_1$. Hence $(u,m+r]\Pi(B,p_i)$ touches $P$ at its
  ending point $p_1'$.
\end{proof}

Theorem~\ref{thm:activities} shows that we can find the internal
and external activities of a basis by just looking at one of its
representations in the diagram $D$. This reduces the problem of
counting the number of bases with given internal and external
activities to the problem of counting the number of lattice paths of a
certain kind in $D$.  In the remainder of this section we sketch a
polynomial-time algorithm that computes this number of bases. 
Note that by Equation~(\ref{tactive}) this yields a different 
approach to computing the Tutte polynomial of a multi-path 
matroid; this approach is slightly quicker than that in the
previous section, but it requires keeping track of more details.

The algorithm uses the characterization of activities in
Theorem~\ref{thm:activities}. Of the conditions in that result,
conditions (b) are somewhat more difficult to deal with; we introduce
the notion of pseudo-activities to count the steps that are active by
conditions (b). Let $R$ be a valid path in the diagram $D$ that ends
in one of the points $p_1',\ldots,p_k'$. Let $s$ be one of its steps
and let $R_s$ be the path that starts at the end of step $s$ and
follows $R$ until its end. We say that $s$ is \emph{pseudo-internally
  active} in $R$ if it is a North step that lies in the top border $Q$
and the path $R_s$ touches the bottom border $P$. Similarly we say
that $s$ is \emph{pseudo-externally active} in $R$ if it is an East
step that lies in $P$ and $R_s$ touches $Q$. Note that, unlike
activities, pseudo-activities are not defined for bases, but for paths
that end at one of the points $p'_1,\ldots,p'_k$ (e.g, the final
segments of paths that correspond to bases).

Let $p$ be a lattice point of the diagram $D$ and let $p_i'$ be one
of the ending points $p_1',\ldots,p_k'$. Let $a$ and $b$ be natural
numbers with $a\leq r$ and $b\leq m$. Let $\tau_P$ and $\tau_Q$ be
variables that can take on the values \emph{true} and \emph{false}. We
define $\Gamma(p,p_i',a,b,\tau_P,\tau_Q)$ to be the number of valid
lattice paths starting at $p$ and ending at $p_i'$ (consisting of one
point if $p=p_i'$), with $a$ pseudo-internally active steps and $b$
pseudo-externally active steps, and touching $P$ if and only if
$\tau_P$ is \emph{true}, and touching $Q$ if and only if $\tau_Q$ is
\emph{true}. The function $\Gamma$ satisfies an easily-verified,
multi-part recurrence relation of which we mention just two parts. Let
$\gamma$ be $\Gamma(p, p_i', a, b, \tau_P, \tau_Q)$ and let $p_E$ and
$p_N$ be, respectively, $p+(1,0)$ and $p+(0,1)$. If $p$ is in neither
$P$ nor $Q$, then
$$\gamma=\Gamma(p_N, p_i', a, b, \tau_P, \tau_Q)+\Gamma(p_E, p_i', a,
b, \tau_P, \tau_Q).$$
If $p$ and $p_N$ are in $Q$, if $p$ is not in
$P$, and if $\tau_Q$ is $true$, then $$\gamma=\Gamma(p_N, p_i',
\bar{a}, b, \tau_P, \tau_Q)+\Gamma(p_E, p_i', a, b, \tau_P,
true)+\Gamma(p_E, p_i', a, b, \tau_P, false)$$
where $\bar{a}$ is
$a-1$ if $\tau_P$ is $true$ and $a$ if $\tau_P$ is $false$, and
$\Gamma(p_N, p_i', \bar{a}, b, \tau_P, \tau_Q)$ is taken to be $0$ if
$\bar{a}<0$ (note, for instance, that this term is also $0$ if
$\bar{a}>0$ and $\tau_P$ is $false$). In this way we get a recurrence
relation that can be expressed in six parts; in each part, $\gamma$ is
a sum of at most three evaluations of $\Gamma$, each involving one of
the points that $p$ leads to, namely, $p_N$ or $p_E$.

With this multi-part recurrence we can compute all values of $\Gamma$
by using a dynamic programming algorithm, not unlike in
Lemma~\ref{lemma:compgraph}. Fix an ending point $p_i'$. Consider a
point $p$ that is $t$ steps away from $p_i'$ and assume we know all
values of $\Gamma$ at 6-tuples involving points that are fewer than
$t$ steps away from $p_i'$. In particular we know all values of
$\Gamma$ at 6-tuples involving $p_E$ and $p_N$, so with the recurrence
relations we can compute any particular value of $\Gamma$ involving
$p$ in constant time. This shows that if we compute the values of
$\Gamma$ for $t$ from $1$ to $r+m$ in this order, then we obtain all
values of $\Gamma$ in $\mathcal{O}(N)$ operations, where $N$ is the
number of 6-tuples in the domain of $\Gamma$, that is, $\mathcal{O}(k
m^2 r^2)$.

We show finally that we can compute the number of bases of internal
activity $i$ and external activity $e$ from $\Gamma$.  This yields a
two-step algorithm for computing the Tutte polynomial of a multi-path
matroid: first compute all values of $\Gamma$, and then obtain the
coefficient of each term $x^iy^e$ in the Tutte polynomial.  The
algorithm requires $\mathcal{O}(km^2r^2)$ operations, or
$\mathcal{O}(n^5)$ where $n$ is $m+r$ (note that $k$ is smaller than
$n$).  This algorithm is somewhat faster than that in
Section~\ref{sec:tutte}.

\begin{lemma}\label{lem:tutte2}
  The number of bases of $M[\ICal]$ with internal activity $i$ and
  external activity $e$ can be found in time $O(k(i+e))$ knowing the
  values of $\Gamma$.
\end{lemma}
\begin{proof}
  We give an algorithm that counts the bases with internal activity
  $i$ and external activity $e$ that contain the element $1$. Note
  that the remaining bases are the complements of the bases of the
  dual with internal activity $e$ and external activity $i$ that
  contain the element $1$, so we can compute their number with the
  same algorithm.
  
  Note that any basis has a unique valid representation that touches
  the top border $Q$, so this gives a one-to-one correspondence between
  bases and certain paths.
  
  For $t>0$ and $j$ in $[k]$, we define $\beta(j,t)$ to be the number
  of bases $B$ such that
\begin{enumerate}
\item[(1)] the internal activity is $i$ and the external activity is
  $e$,
\item[(2)] $[t]\subseteq B$ and $t+1\notin B$, and
\item[(3)] the path $\Pi(B,p_j)$ is the unique valid representation
  that touches $Q$.
\end{enumerate}

Let $R$ be the path $N^t E$ that starts at $p_j$ and let $s$ be the
last step of $R$. By conditions (2) and (3), $R$ coincides with
$[1,t+1]\Pi(B,p_j)$ where $B$ is any of the bases that $\beta(j,t)$ is
counting; clearly if $R$ is not valid, then $\beta(j,t)$ is $0$.  The
first $t$ North steps of $R$ are internally active elements in $B$,
but the step $s$ may or may not be externally active. If $s$ does not
lie in the bottom border $P$, then by Theorem~\ref{thm:activities} it
is not externally active, so
$$\beta(j,t)=\sum_{ (\tau_P,\tau_Q)\in T_P\times T_Q } \Gamma(p, p_j',
i-t, e, \tau_P, \tau_Q),$$
where $p$ is the point $p_j+(1,t)$ where
$R$ ends, $T_P$ is $\{true, false\}$ and $T_Q$ is $\{true\}$ if $R$
does not touch $Q$, and $\{true, false\}$ if $R$ touches it.  The two
possibilities for $T_Q$ arise from the requirement in condition (3)
that the paths touch $Q$, so if $R$ does not touch it, then the
remaining part of the path has to. Notice that many terms in this sum
may be $0$; for instance, if $i-t$ or $e$ are greater than $0$, then
every evaluation of $\Gamma$ where $\tau_P$ or $\tau_Q$ is $false$ is
$0$.

Now assume that $s$ lies in $P$. (This can happen only if $j$ is $1$.)
Note that $s$ is externally active if and only if the remaining part
of the path touches $Q$.  Therefore 
$$\beta(j,t)=\sum_{ (\tau_P,\tau_Q)\in T_P\times T_Q } \Gamma(p, p_j',
i-t, e-\delta(\tau_Q), \tau_P, \tau_Q),$$
where $p$ is the point
$p_j+(1,t)$, the set $T_Q$ is $\{true,false\}$ or $\{true\}$ according
to whether $R$ does or does not touch $Q$, the set $T_P$ is
$\{true\}$, and $\delta(\tau_Q)$ is $1$ if $\tau_Q$ is $true$, and $0$
otherwise. Note that since $s$ lies in $P$ the paths we are counting
start on $P$ so $\tau_P$ must be $true$.

To obtain the number of bases with internal activity $i$ and
external activity $e$ containing $1$ we add up all the values of
$\beta$. Since the number $t$ is bounded by $i$, we can do the
computations in $O(ki)$ operations. The same algorithm applied to
the dual matroid needs $O(ke)$ operations, hence we can compute
the total number of bases of $M[\ICal]$ of internal activity $i$
and external activity $e$ from $\Gamma$ in $O(k(i+e))$ operations.
\end{proof}

\section{Further Structural Properties}\label{sec:structure}

This final section treats a variety of properties of multi-path
matroids and their presentations.

Every connected lattice path matroid with at least two elements has a
spanning circuit~\cite[Theorem~3.3]{lpm2}.  The analogous property
holds for multi-path matroids, as we now show.  By the result just
cited, it suffices to focus on multi-path matroids that are not
lattice path matroids.

\begin{thm}\label{thm:spct}
  A multi-path matroid $M[\mathcal{I}]$ that is not a lattice path
  matroid has a spanning circuit.  Furthermore, every element is in
  some spanning circuit.
\end{thm}

\begin{proof}
  Since multi-path matroids of rank less than $2$ are lattice path
  matroids, we are assuming that the rank is at least $2$.  The set
  $F=\{f_I\,:\,I\in\mathcal{I}\}$ of first elements is a proper subset
  of the ground set $S$.  By the comments at the end of
  Section~\ref{sec:def}, the first element $f_I$ of any interval $I$
  in $\mathcal{I}$ is in both $I$ and $\Sigma^{-1}(I)$; also,
  $M[\mathcal{I}]$ has no loops.  From these observations, it is
  immediate to check that $F\cup x$, for any $x$ in $S-F$, is a
  spanning circuit of $M[\mathcal{I}]$.
\end{proof}

Corollary~\ref{cor:con} follows from Theorem~\ref{thm:spct} since
multi-path matroids that are not lattice path matroids have no loops,
and loopless matroids with spanning circuits are connected.

\begin{cor}\label{cor:con}
  Every multi-path matroid that is not a lattice path matroid is
  connected.
\end{cor}

From Corollary~\ref{cor:con}, or directly from Theorem~\ref{thm:spct},
it follows that, in contrast to the class of lattice path matroids,
the class of multi-path matroids is not closed under direct sums.  For
example, recall that the $3$-whirl $\mathcal{W}^3$ is a multi-path
matroid but not a lattice path matroid.  Therefore the direct sum
$\mathcal{W}^3\oplus \mathcal{W}^3$ is neither a lattice path matroid
(since $\mathcal{W}^3$ is a restriction) nor a multi-path matroid.  Of
course, one could consider the class of matroids whose connected
components are multi-path matroids; such matroids can be realized with
a simple variation on Definition~\ref{def:multipath}, having the
intervals being intervals in the cycles in the cycle decomposition of
an arbitrary permutation of the ground set.

The next theorem gives some indication of how close multi-path
matroids are to lattice path matroids.

\begin{thm}\label{thm:flat}
  The restriction of a multi-path matroid to a proper flat is a
  lattice path matroid.
\end{thm}

\begin{proof}
  Let $M$ be the multi-path matroid $M[\mathcal{I}]$.  The class of
  lattice path matroids is closed under direct
  sums~\cite[Theorem~3.6]{lpm1}, so it suffices to prove the assertion
  for proper flats $F$ for which $M|F$ is connected.  The assertion is
  easily seen to hold for flats of rank $2$ or less.  Let $F$ be a
  proper flat of rank $3$ or more for which $M|F$ is connected.  By
  Theorem~\ref{thm:spct} and the corresponding result for lattice path
  matroids~\cite[Theorem~3.3]{lpm2}, the restriction $M|F$ has a
  spanning circuit $C$.  It follows from Hall's matching theorem that
  a circuit $C'$ of a transversal matroid has nonempty intersection
  with exactly $|C'|-1$ of the sets in any presentation; therefore the
  inequality $|C|-1=r(F)<r(M)$ implies that $C$ is disjoint from at
  least one interval $I$ in $\mathcal{I}$.  Thus, $F$ is disjoint from
  $I$, so $F$ is a flat of the deletion $M\del I$.  Observe that
  $M\del I$ is a lattice path matroid: by Lemma~\ref{lem:trimpres},
  the presentation $(J\del I\,:\,J\in \mathcal{I},J\ne I)$ of $M\del
  I$ by intervals in
  $\sigma(l_I),\sigma^2(l_I),\ldots,\sigma^{-1}(f_I)$ contains a
  presentation of $M\del I$ by an antichain of intervals.  Since
  $M\del I$ is a lattice path matroid, so is $M|F$.
\end{proof}

Theorem~\ref{thm:flat} allows one to carry over certain results about
lattice path matroids to multi-path matroids.  For instance, the
description of the circuits of lattice path
matroids~\cite[Theorem~3.9]{lpm2} applies to the nonspanning circuits
of multi-path matroids.  We mention several other results that are
counterparts of results for lattice path matroids and that may prove
useful for the further study of multi-path matroids.  Let
$M[\mathcal{I}]$ be a multi-path matroid of rank $r$ on the set $S$.
\begin{enumerate}
\item Let $I_{i_1},I_{i_2},\ldots,I_{i_h}$ be the intervals in
  $\mathcal{I}$ that have nonempty intersection with a fixed connected
  flat $F$ of $M[\mathcal{I}]$ of rank greater than $1$.  Then
  $\{I_{i_1},I_{i_2},\ldots,I_{i_h}\}$ is a $\Sigma$-interval in
  $\mathcal{I}$ and $h$ is $r(F)$.
\item Statement (1) implies that there are at most $r$ connected flats
  of a fixed rank greater than $1$ in $M[\mathcal{I}]$. Whirls show
  that this bound cannot be improved.
\item Statement (1) also implies that any flat of $M[\mathcal{I}]$ is
  covered by at most two connected flats.
\item The elements in any connected flat of $M[\mathcal{I}]$ form a
  $\sigma$-interval in $S$.
\item If $X_1,X_2,$ and $X_3$ are connected flats of $M[\mathcal{I}]$,
  and if no two sets among $X_1,X_2,X_3$ are disjoint, then either one
  of $X_1,X_2,X_3$ is contained in the union of the other two, or
  $X_1\cup X_2\cup X_3$ is $S$.
\end{enumerate}
(Compare statements (1) and (4) with~\cite[Theorem~3.11]{lpm2};
compare statements (2) and (3) with~\cite[Corollary~3.12]{lpm2}.)

Our final topic is minimal presentations of multi-path matroids.
Example 3 in Section~\ref{sec:def} gives distinct $\sigma$-interval
presentations of a multi-path matroid that are also minimal
presentations.  The next theorem shows that any minimal
$\sigma$-interval presentation is also a minimal presentation.  Note
that the converse is not true: for example, the presentation
$(\{1,4\},\{2,4\},\{3,4\})$ of $U_{3,4}$ is minimal but these sets are
not $\sigma$-intervals for any cycle $\sigma$ on $[4]$.

\begin{thm}\label{thm:min}
  The sets in a minimal $\sigma$-interval presentation of a multi-path
  matroid are cocircuits of the matroid.  Any minimal
  $\sigma$-interval presentation of a multi-path matroid is a minimal
  presentation.
\end{thm}

\begin{figure}
\begin{center}
\includegraphics[width = 3.8truein]{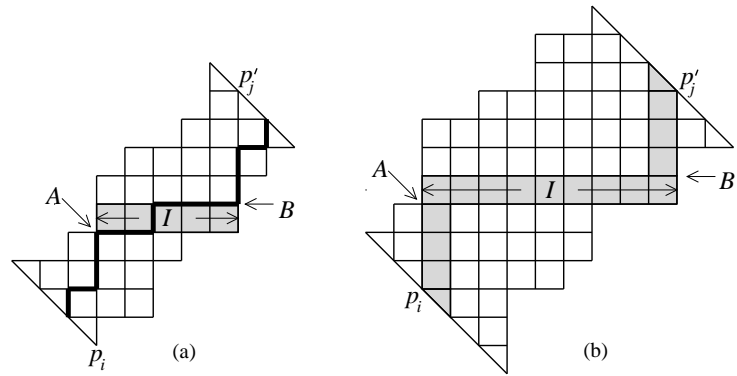}
\end{center}
\caption{The cases (a) $i\leq j$ and (b) $i>j$ in the proof of
  Theorem~\ref{thm:min}. }\label{minimal}
\end{figure}

\begin{proof}
  Assume that the multi-path matroid $M$ has rank $r$ and that
  $\mathcal{I}$ is a minimal $\sigma$-interval presentation of $M$.
  Each set in a presentation of a transversal matroid is the
  complement of a flat of the matroid.  Since cocircuits are the least
  nonempty complements of flats, a presentation by cocircuits is
  necessarily minimal, so the second assertion of the theorem follows
  from the first.  Let $I$ be in $\mathcal{I}$.  Since the complement
  of $I$ is a flat, the first assertion follows if we show that this
  complement contains $r-1$ independent elements.  In terms of lattice
  paths, we need to show that there is a lattice path in some diagram
  $D(\mathcal{I},x)$ that connects a pair of corresponding points
  $p_h$ and $p'_h$ and has only one North step that is labelled by an
  element of $I$.  This statement is trivial if $r$ is $1$, so assume
  $r$ exceeds $1$.
  
  Since $\mathcal{I}$ is an antichain and $r$ exceeds $1$, some
  element, say $x$, of $M$ is not in $I$.  Let $A$ and $B$ be,
  respectively, the lower left and upper right points in the row of
  $D(\mathcal{I},x)$ that represents $I$.  (See Figure~\ref{minimal}.)
  Let $i$ be the least positive integer for which there is a path in
  $D(\mathcal{I},x)$ from $p_i$ to $A$.  Note that there is a path in
  $D(\mathcal{I},x)$ from $p_h$ to $A$ if and only if $h\geq i$.
  Similarly, let $j$ be the greatest integer for which there is a path
  from $B$ to $p'_j$.  Thus, there is a path in $D(\mathcal{I},x)$
  from $B$ to $p'_h$ if and only if $h\leq j$.  It follows that if
  $i\leq j$, then there is a path in $D(\mathcal{I},x)$ that connects
  any pair of corresponding points $p_h$ and $p'_h$ with $i\leq h\leq
  j$ and that has only one North step labelled by an element of $I$,
  as desired.  We complete the proof by showing that the alternative,
  the inequality $i>j$, contradicts the assumption that $\mathcal{I}$
  is a minimal $\sigma$-interval presentation.  The inequality $i>j$
  forces $i$ to be greater than $1$.  If there were a path in
  $D(\mathcal{I},x)$ of the form $N^aEQ$ from $p_i$ to $A$, then the
  path $N^{a+1}Q$ from $p_{i-1}$ to $A$ would also be in
  $D(\mathcal{I},x)$, contrary the choice of $i$, so there is only one
  path from $p_i$ to $A$ and this path consists of all North steps.
  Similarly, $j<k$ and the unique path from $B$ to $p'_j$ consists of
  all North steps.  From these conclusions, it follows that for any
  path in $D(\mathcal{I},x)$, say from $p_h$ to $p'_h$, that uses the
  North step labelled $f_I$ in the row corresponding to $I$, or any
  North step immediately above this one, we have $h\geq i$ and the
  same sequence of steps, but instead going from $p_{h-1}$ to
  $p'_{h-1}$, remains in $D(\mathcal{I},x)$.  Thus, by deleting $f_I$
  from $I$, deleting $\sigma(f_I)$ from $\Sigma(I)$ if
  $f_{\Sigma(I)}=\sigma(f_I)$, deleting $\sigma^2(f_I)$ from
  $\Sigma^2(I)$ if $f_{\Sigma^2(I)}=\sigma^2(f_I)$, etc., we obtain a
  smaller $\sigma$-interval presentation of $M$, that, as desired,
  contradicts the assumed minimality of $\mathcal{I}$.
\end{proof}

Let $M$ be a matroid of rank $r$ and nullity $m$.  Since any
hyperplane contains at least $r-1$ of the $r+m$ elements of $M$, any
cocircuit has at most $m+1$ elements.  From this observation, the
following corollary of Theorem~\ref{thm:min} is evident.

\begin{cor}
The sets in any minimal $\sigma$-interval presentation of a multi-path
matroid of nullity $m$ have at most $m+1$ elements.
\end{cor}

\begin{center}
\textsc{Acknowledgements}
\end{center}

The authors thank Anna de Mier and Marc Noy for useful discussions
about the material in this paper and the exposition.

\end{document}